\documentclass{article}[11pt]                     

\usepackage{geometry}                
\usepackage{a4wide}
\usepackage{graphicx}
\usepackage{amssymb}
\usepackage{color,paralist,amsmath}

\newcommand{\C}{\mathcal{C}}
\newcommand{\G}{\mathcal{G}}

\renewcommand{\P}{\mathcal{P}}
\renewcommand{\S}{\mathcal{S}}
\newcommand{\U}{\mathcal{U}}

\newtheorem{theorem}{Theorem}[section]

\newtheorem{definition}[theorem]{Definition}

\begin{document}

\title{Design Heuristic for Parallel Many Server Systems}


\author{Ivo Adan\thanks{
Department of Industrial Engineering, Eindhoven University of Technology, P.O. Box 513, 5600 MB Eindhoven, the Netherlands; email iadan@tue.nl
} \and Marko Boon\thanks{
Department of Mathematics and Computer Science, Eindhoven University of Technology, P.O. Box 513, 5600 MB Eindhoven, the Netherlands; email m.a.a.boon@tue.nl
}
\and  Gideon Weiss\thanks{
Department of Statistics,
The University of Haifa,
Mount Carmel 31905, Israel; email
gweiss@stat.haifa.ac.il
Research supported in part by
Israel Science Foundation Grants 711/09 and 286/13.}
}

\date{\today}
\maketitle
\begin{abstract}
We study a parallel queueing system with multiple types of servers and customers.
A bipartite graph describes which pairs of customer-server types are compatible. We consider the
service policy that always assigns servers to the first, longest waiting compatible customer, and that always assigns customers to the longest idle compatible server if on arrival, multiple compatible servers are available.
For a general renewal stream of arriving customers and general service time distributions, the behavior of such systems is very complicated. In particular, the calculation of matching rates, the fraction of services of customer-server type, is intractable. We suggest through a heuristic argument that if the number of servers becomes large, the matching rates are well approximated by matching rates calculated from the tractable bipartite infinite matching model. We present simulation evidence to support this heuristic argument, and show how this can be used to design systems with desired performance requirements.

\end{abstract}

\section{Introduction}
\label{sec.introduction}
Parallel service systems have servers of types $\S=\{s_1,\ldots,s_J\}$, customers of types
$\C=\{c_1,\ldots,c_I\}$, and bipartite compatibility graph $\G\subseteq \C \times \S$, where $(c_i,s_j)\in \G$ if servers of type $s_j$ can serve customers of type $c_i$.  They model situations in which a large volume of service requests of various types are channelled to a central facility, where they are attended by a large number of agents differentiated by skill.  Such situations commonly occur in manufacturing, transportation, service contact centers, health systems, communications, internet data exchange, computing and various other areas of applications.
A queueing model for this has a general renewal stream of arriving customers with rate $\lambda$, where successive arrivals are of i.i.d types, $c_i$ with probability $\alpha_{c_i}$, and there is a total of $n$ servers, $n_{s_j}$ of which are of type $s_j$.  Service times are independent, distributed according to general distributions $G_{c_i,s_j}$, with mean $m_{c_i,s_j}$ and service rate $\mu_{s_j,c_i}=1/m_{c_i,s_j}$.  Customers have finite patience, with independent patience time distributions $F_{c_i}$, and a customer abandons if he does not start service  by the time his patience  is exhausted.

Parallel server systems are widely discussed in the literature.  An incomplete list would include an early study \cite{green:85};  applications to manufacturing and supply chain management    \cite{veeger-etman-rooda:08,rubino-ata:09},  applications to call centers and internet service systems \cite{gans-koole-mandelbaum:03,harchol-balter-etal:99,squillante-etal:01,wallace-whitt:05},  attempts to find optimal policies, mainly for small graph systems \cite{williams:00,bell-williams:01,armony-ward:10,armony-ward:13,ghamami-ward:13,tezcan-dai:10}, heavy traffic and fluid approximations \cite{harrison-lopez:99,harrison-zeevi:05}, and many server scaling \cite{gurvich-whitt:09,gurvich-whitt:10}.
Most relevant to our current paper are \cite{foss-chernova:98,adan-weiss:14,talreja-whitt:07,nov-weiss-zhang:16}.

In assessing such systems there are various objectives that may be of importance, on the customer side they include waiting times and abandonment rates as well as consideration of  fairness to customers of various types or priorities for some types.  In conflict with those, on the server side there is the objective of maximum utilization of the servers, minimizing their number, and reaching a balanced work division between the various types.   Each of these may carry a different weight in different application contexts.   Often $\lambda$, $\alpha_{c_i}$ and $F_{c_i}$ are given, together with some form of  quality of service requirements.  All the other parameters of the system can be adjusted to achieve the requirements in an optimal way:  One can redesign the bipartite compatibility graph, change the service rates, change the workforce mix, decide on $n$, and decide on the service policy.  It should perhaps be pointed out that changing the service policy may be as hard and costly as adjusting any of the other service parameters.   At this level of generality such systems do not allow a complete analytic analysis, and performance is often evaluated in practice by simulation.  However, any methods for calculating approximate performance measures or supporting design without the need to use simulation should be quite valuable. It is the aim of this paper to deliver such methods.

In the current paper we focus on the policy of first come first served (FCFS), where whenever a server is available he will take the longest waiting compatible customer, and assign longest idle server (ALIS), where whenever a customer arrives he will be assigned to the longest idling compatible server.  We provide a heuristic to calculate performance measures under this policy, when $\lambda$ and $n$ are large.

FCFS-ALIS in a parallel service system has several advantages:   It attempts to achieve resource pooling \cite{tsitsiklis-xu:12}, i.e. all the servers are busy for about the same fraction of time, and it attempts to give all customers the same service level, i.e. global FCFS, equally for all types of customers \cite{talreja-whitt:07}.  It is also fair to the servers.  One notable property of FCFS is the following:
Assume that arriving customers can choose the server they wish to go to, and each server then serves his queue FCFS.  If each arrival has complete information on the schedule of all the servers at his moment of arrival, then to minimize his waiting time he will join the
compatible server that has the shortest workload (JSW).  Thus JSW is the Nash equilibrium of fully informed customers minimizing waiting times.  But this policy of JSW is automatically achieved when customers queue up in a single queue and the servers are using FCFS.
FCFS can then serve as a benchmark, and comparison of the costs under FCFS with other policies will provide an estimate of the price of anarchy.
Apart form that, FCFS is easy to implement, as it does not require any online calculations or knowledge of system parameters.  It is also sometimes required by law.  Finally it is indeed a policy very commonly used in practice.  On the negative side, FCFS may waste resources by letting servers serve customers for which they are not efficient, and it may cause long delays to customer types that have a limited number of compatible servers.  However, some of these shortcomings can be avoided by redesigning the compatibility graph.  It is safe to say that in practice most service systems use FCFS for a large proportion of their operation, even if they implement some more sophisticated policies in some of their service decisions.

Unfortunately, analysis of parallel service systems under FCFS is very hard.  Foss and Chernova \cite{foss-chernova:98} provide an example of a symmetric system with 3 types of customers and 3 servers, and just 2 service distributions with fixed fast and slow service rates, where stability of the system depends on the entire shape of the service time distributions.  The difficulty is in calculating the matching rates $r_{c_i,s_j}$, defined as the long term average fraction of customers of type $c_i$ which are processed by servers of types $s_j$.  Given the matching rates, one can calculate the total service capacity of the system, as
\[
\mu = \sum_{(c_i,s_j)\in \G}  r_{c_i,s_j} \mu_{c_i,s_j},
\]
and then conclude that under FCFS the system is stable if $\lambda<\mu$.
It is the calculation of the matching rates, how many customers of type $c_i$ are served by servers of type $s_j$ under FCFS-ALIS, which is intractable, and may depend on the entire shape of the service time distributions.   Matching rates can be calculated for some types of graphs \cite{talreja-whitt:07,nov-weiss-zhang:16}, and they can also be calculated for general bipartite graphs when arrivals are Poisson, service rates depend only on the servers, and services are exponential \cite{adan-weiss:14}, but not otherwise.  However, matching rates can be calculated for the much simpler and very tractable FCFS infinite bipartite matching  model \cite{caldentey-kaplan-weiss:09,adan-weiss:11,adan-busic-mairesse-weiss:15}.

It is our conjecture that  when $\lambda$ and $n$ are large, the matching rates of the general parallel service system under FCFS-ALIS are approximated by those of the FCFS infinite bipartite matching model: this is the basis for our current paper. A first step towards verifying the conjectured many server behavior is made in \cite{zhan-weiss:16} by studying the limiting behavior of the ``N'' system. Even for this very simple system, the derivations of the limits are quite laborious, emphasizing the difficulty in verifying the conjecture for the general system. We use our ability to calculate matching rates in order to design parallel service systems operating under FCFS-ALIS.
In \cite{adan-boon-busic-mairesse-weiss:13,adan-boon-weiss:13} we presented a heuristic algorithm to determine the required workforce in the Efficiency Driven (ED) overloaded regime. The main contribution of the current paper is to extend this algorithm to other regimes of interest: The Quality Driven (QD) mode in which there is underload, and Quality and Efficiency Driven (QED) mode in which the load is exactly one \cite{mandelbaum-zeltyn:05}.
So we will consider parallel service systems operating in ED mode, QD mode and QED mode. Our objective in each of these is to design the workforce required to achieve certain service requirements, specifically:
 \begin{compactitem}
\item[-]
In ED mode, we specify average waiting times for customers, and resulting abandonment rates.
\item[-]
In QD mode, we specify average idle time for the servers.
\item[-]
In QED mode we design the workforce to achieve almost full utilization, zero or short waiting times, and no abandonments.
\end{compactitem}
Under FCFS-ALIS we achieve these pre-specified requirements with complete resource pooling of servers, and balanced service levels for all types of customers.

We also present designs where under FCFS-ALIS the parallel servers are not pooled, and use this to achieve differentiated service levels for the various types of customers, based on pre-specified priority levels.  In this case the design achieves the desired service levels and utilizations for all types, under the policy of FCFS-ALIS, without the need to use prioritized service decisions.

In deciding on the staffing levels we consider two choices:  We can specify the fraction of services provided by each type of server, out of all the services, we denote this by $\beta_{s_j}$.  This specification is appropriate if for example
costs are calculated per service, and may differ for different types of servers.  Alternatively we may specify the fraction of total number of servers of each type, we denote this by $\theta_{s_j}=n_{s_j}/n$.  This specification is appropriate if for example the costs are per server, and may differ for different types of servers.
We consider staffing decisions made on the basis of either of these choices.

We stress that our purpose in this design is not to minimize the number of servers, or to minimize  waiting times for a specified number of servers.
Such minimization  will usually be achieved by a tree compatibility graph, for which matching rates are easily computed (cf. \cite{nov-weiss-zhang:16}).
Rather it is to start from a given compatibility graph, and  obtain designs which achieve resource pooling and balance of service.  The advantage of general  compatibility graphs is an added robustness against variable arrival rates and customer type composition, and an added flexibility in the quality of service.

The rest of the paper is structured as follows.  In Section \ref{sec.infinitematching} we describe the FCFS infinite bipartite matching model, and the formula for the calculation of matching rates.  In Section \ref{sec.manyserver} we present our conjecture on the behavior of FCFS-ALIS parallel service systems under many server scaling, that links them with the infinite bipartite matching model.  In Section \ref{sec.algorithm}  we present our design algorithms, based on the calculation of matching rates.   In Section \ref{sec.examples} we present examples in which we calculate designs, and examine the performance under our designs.
We present extensive simulation results, that confirm the validity of our approach for a range of $\lambda$ and $n$ scales.  Finally, in Section \ref{sec.discussion} we discuss significance of our results and how we think they should be applied in practice.


\section{FCFS  infinite bipartite matching}
\label{sec.infinitematching}

\begin{figure}[htb]
   \begin{center}
    \includegraphics[scale=0.35]{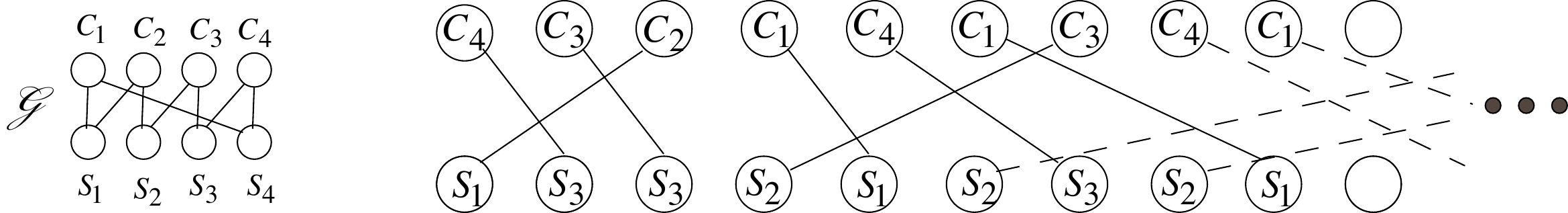}
    \caption{Matching process for given $4 \times 4$ graph $\G$.}\label{fig0}
   \end{center}
\end{figure}

We now consider a system with customer types $\C$ and server types $\S$, with a bipartite compatibility graph $\G$, and a much simplified stochastic model:  We have  infinite sequences of customers $c^1,c^2,\ldots,c^m,\ldots$ where $c^m\in \C$ and of servers $s^1,s^2,\ldots,s^n,\ldots$ where $s^n\in \S$.  We assume that $c^m$ are drawn  according to probabilities $\alpha=(\alpha_{c_1},\ldots,\alpha_{c_I})$ and $s^n$ are drawn according to probabilities $\beta=(\beta_{s_1},\ldots,\beta_{s_J})$, and they are all independent.   For each realization of the sequences we match customers and servers according to a FCFS
policy:  $s^n$ is matched to the earliest compatible $c^m$ in $c^1,c^2,\ldots,$ which has not yet been matched to $s^1,\ldots,s^{n-1}$. The matching process, for given graph $\G$, is illustrated in Figure \ref{fig0}.
This model  is much simpler than a queueing model, since it involves no arrival times, no service times, no busy or idle servers, and since it treats customers and servers in an entirely symmetric way.  This system is studied in \cite{caldentey-kaplan-weiss:09,adan-weiss:11,adan-busic-mairesse-weiss:15}.  It is shown in \cite{adan-weiss:11} that the matching is uniquely determined for any two sequences and that all customers and servers are matched almost surely.  Furthermore, the system demonstrates dynamic reversibility, and is associated with a Markov chain that has a product form stationary distribution.  The stationary distribution is used to obtain explicit expressions for the matching rates.  We describe the calculation of the matching rates now.

We use the following notations:  we let $\C(s_j)$ be the set of customer types compatible with server type $s_j$, and $\S(c_i)$ be the set of server types compatible with customer type $c_i$.  For a subset of customer types $C$ we let $\S(C)=\bigcup_{c_i\in C} \S(c_i)$, and for a subset of server types $S$ we let $\C(S)=\bigcup_{s_j\in S} \C(s_j)$.  We also let $\U(S)=\overline{\C(\overline{S})}$ be the customer types that can only be served by servers of types in $S$.
For subsets $C,\,S$ we define $\alpha_C=\sum_{c_i\in C} \alpha_{c_i}$, and $\beta_S = \sum_{s_j\in S} \beta_{s_j}$.
%
%
%
\begin{definition}
For given $\alpha,\beta,\G$ we say that there is complete resource pooling in the FCFS infinite bipartite matching system if the following three equivalent conditions hold:
\begin{equation}\label{eqn.crp}
 \alpha_C < \beta_{\S(C)},  \quad  \beta_S < \alpha_{\C(S)}, \quad  \beta_S > \alpha_{\U(S)}, \qquad
 S\subset \S,\, S\ne\emptyset,\S,  \quad C\subset \C,\,C\ne \emptyset,\C.
\end{equation}
\end{definition}
\begin{theorem}[from \cite{adan-weiss:11}]
Let $r_{c_i,s_j}(n)$ be the (random) number of $c_i,s_j$ matches  between $c^1,\ldots,c^n$ and $s^1,\ldots,s^n$, in the FCFS infinite bipartite matching of the two sequences.
If complete resource pooling holds, then almost surely $\lim_{n\to\infty} r_{c_i,s_j}(n) = r_{c_i,s_j}$ which is calculated by
\begin{eqnarray}
\label{eqn.matches}
r_{c_i,s_j} &=&  \beta_{s_j} \sum_{\mathcal{P}_J}  B \prod_{k=1}^{J-1} (\beta_{(k)} - \alpha_{(k)})^{-1}   \nonumber \\
&& \left(
\sum_{k=1}^{J-1} \phi_k
\frac{\alpha_{(k)}}{\beta_{(k)}-\alpha_{(k)} \chi_{k}}
\prod_{l=1}^{k-1} \frac{\beta_{(l)}-\alpha_{(l)}}{\beta_{(l)}-\alpha_{(l)} \chi_{l}}
+ \frac{\phi_J}{\phi_J+\psi_J}
\prod_{l=1}^{J-1} \frac{\beta_{(l)}-\alpha_{(l)}}{\beta_{(l)}-\alpha_{(l)} \chi_{l}}
\right) ,
\end{eqnarray}
where the summation is over $\P_J$, the set of all permutations of the server types  $\S$, and for each permutation of the servers $S_1,\ldots,S_J$, the following notation is used:
\[
\alpha_{(k)} = \alpha_{U\{S_1,\ldots,S_k\}}, \qquad \beta_{(k)} =
\beta_{\{S_1,\ldots,S_k\}},\qquad k=1,\ldots,J,
\]
\[
\phi_k = \frac{\alpha_{U\{S_1,\ldots,S_k\}\cap \{c_i\}}}{\alpha_{U\{S_1,\ldots,S_k\}}}, \qquad
\psi_k = \frac{\alpha_{U\{S_1,\ldots,S_k\}\cap (C(s_j)\backslash\{c_i\})}}{\alpha_{U\{S_1,\ldots,S_k\}}}, \qquad
\chi_k=1-\phi_k-\psi_k \ ,
\]
and $B$ is the normalizing constant:
\[
B^{-1} =\sum_{\mathcal{P}_J}
\left(
 (\beta_{\{S_1\}} -  \alpha_{\U\{S_1\}}) (\beta_{\{S_1,S_2\}} -  \alpha_{\U\{S_1,S_2\}})
\cdots (\beta_{\{S_1,\ldots,S_{J-1}\}} -  \alpha_{\U\{S_1,\ldots,S_{J-1}\}})
\right)^{-1} .
\]
\end{theorem}

An easy example of this formula is for the case that $I=J$ and $\C(s_j)=\C\backslash c_j$, i.e the bipartite compatibility graph is almost complete, each server can serve all but one of the customer types.  In that case, complete resource pooling holds if and only if $\alpha_{c_j}+\beta_{s_j} < 1$, and the matching rates are:
\begin{eqnarray}
\label{eqn.almostcomplete}
r_{c_i,s_j} &=&
\alpha_{c_i} \beta_{s_j}
\frac{
(1- \alpha_{c_i})(1- \beta_{s_j}) - \alpha_{c_j} \beta_{s_i}}
{(1- \alpha_{c_i} - \beta_{s_i})(1- \alpha_{c_j} - \beta_{s_j})}  \Big{/}
 \left(1+\sum_{i=1}^I {\alpha_i \, \beta_i  \over
1- \alpha_i - \beta_i }\right).
\end{eqnarray}
However, for any other  bipartite compatibility graph, the formula (\ref{eqn.matches}) does not seem to simplify, and we suspect that its calculation is $\sharp P$ hard.  We have programmed it to be able to calculate it up to $I,J \le 12$, but it will become hard to compute the matching rates for larger number of types.
Recently, \cite{} developed a highly accurate and efficient approximation for the matching
rates $r_{c_i,s_j}$ based on Ohm’s Law (which in some cases reduces to exact results).
This provides an attractive alternative to the exact solution in case of many customer and server types.  .

When resource pooling does not hold, it is shown in \cite{adan-weiss:14} that there is a unique decomposition $(\C,\S)$ into subsystems $(\C^{(1)},\S^{(1)}),\ldots,(\C^{(L)},\S^{(L)})$, such that
\begin{equation}
\label{eqn.decomposition}
\frac{\beta_{\S^{(1)}}}{\alpha_{\C^{(1)}}} < \cdots < \frac{\beta_{\S^{(L)}}}{\alpha_{\C^{(L)}}},  \qquad
(\C^{(l)},\S^{(l)}) \mbox{ has complete resource pooling, } l=1,\ldots,L.
\end{equation}

A Mathematica program to exactly calculate the matching rates for given $\alpha,\beta,\G$ is available from the authors.


\section{Matching under many server scaling}
\label{sec.manyserver}
Consider a queueing system with a single customer type  and a single server type, with arrival rate $\lambda$, and patience distribution $F$,  and with $n$ servers, each with service rate $\mu$, so that the traffic intensity is $\rho = \lambda/n\mu$.  Many server scaling occurs when we keep $\mu$ and $\rho$ fixed and let both $\lambda$ and $n$ increase. Note that, to increase $\lambda$, we scale the inter-arrival time distribution, and thus we do not alter its shape. Because of abandonments the system will always be stable.  There will be three behavior modes for this system:  When $\rho<1$ the system is in QD ({\em quality driven mode}).  In QD mode, there is always a fraction $\approx\! (1-\rho)$ of idle servers and customers never wait and nobody abandons.  When $\rho>1$ the system is in ED ({\em efficiency driven mode}).  In ED mode, servers are always busy,  there is always a queue, and a fraction $\approx\! F(W) = (\rho - 1)/\rho$ of customers abandons without service.    Customers  with patience $\le W$ do not get served, and customers with patience $> W$ receive service after a wait of $\approx\! W$.  When $\rho\approx\! 1$ the system is in QED ({\em quality and efficiency driven mode}).  In QED mode, servers are busy most of the time and if they idle it is only for a short while, an appreciable  fraction of customers do not need to wait, the remaining customers wait a very short time, and very few customers abandon~\cite{whitt:06}.

We now consider the system of parallel skill based servers of Section \ref{sec.introduction}.  We fix the fractions $\alpha_{c_i}$, $\theta_{s_j} = n_{s_j}/n$, the service time distributions $G_{c_i,s_j}$, and the patience distributions $F_{c_i}$.  We use FCFS-ALIS policy, and we let $\lambda$ and $n$ increase at the same rate, so that we get into many server scaling.   We cannot directly calculate $\rho$ for this system, as it depends on the service policy, and in particular we cannot calculate it directly under FCFS-ALIS, since it depends on the matching rates, but we will try and approximate it.  Under many server scaling we can expect that for favorable choices of parameters, the system will achieve resource pooling, so that customers of different types will have similar waiting times,  and servers of different types will have similar workloads and similar idle times.  Under such conditions, the system will again behave in one of the three modes, QD, ED or QED, according to the traffic intensity.

We now make the following conjectures regarding the behavior of the system under many server scaling, when complete resource pooling holds:  First, we conjecture that the order in which customers will reach the head of the line (if they did not abandon previously) will be such that the types of customers will be i.i.d. with some probabilities $\tilde{\alpha}_{c_i}$, approximately.  Next,  we conjecture that the order in which servers will become available at completion of service, will be such that the types of servers will be i.i.d with some probabilities $\beta_{s_j}$, approximately.   The rationale  for the first conjecture is that under complete resource pooling, when $\lambda$ and $n$ are large,
under QD or QED customers will not wait at all and so the successive types will be i.i.d with probabilities $\alpha_{c_i}$, while under ED,
 all customers will need to wait for approximately a constant time $W$ before being served, and so a fraction $F_{c_i}(W)$ of customers of type $c_i$ will abandon,
 and successive customer types that reach the head of the line will be i.i.d. with some modified probabilities.
The rationale  for the second conjecture is that in steady state the time points at which each individual server will complete a job form a stationary point process, and with many servers these point processes will be nearly independent.  It is known that the superposition of many independent stationary point processes, under appropriate scaling, converges to a Poisson process (cf. Khinchine \cite{khinchine:60}), which indicates that our conjecture on server availability may in fact be true.
Under these two conjectures, the matching between customers and servers will be approximately the same as for the FCFS infinite bipartite matching model.

The following three figures illustrate the operation of our system under FCFS-ALIS in each of the above mentioned three modes:

\begin{figure}[htb]
  \begin{center}
    \includegraphics[scale=0.40]{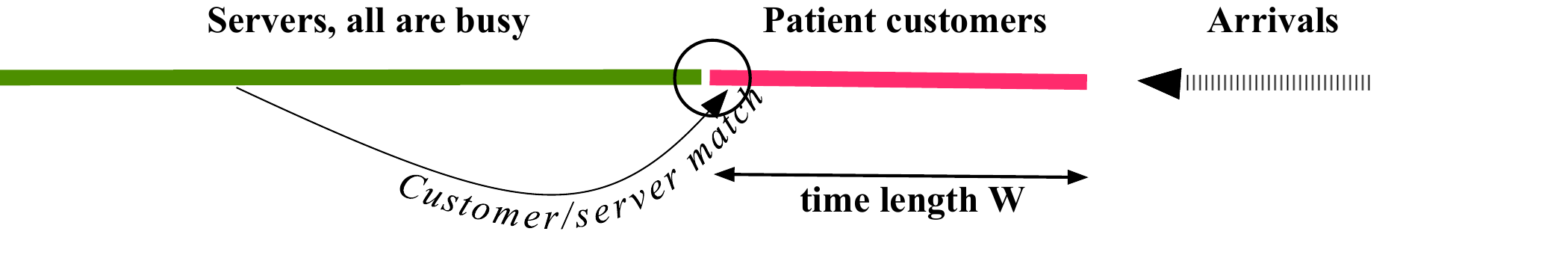}
   \end{center}
   \caption{Operation of system in ED mode}\label{fg1}
\end{figure}

In ED mode, all the servers are always busy, customers with enough patience wait a time $W$, and when they reach the head of the queue, they match with the next compatible server (see Figure \ref{fg1}).
Note that customers entering service are still of i.i.d. types, approximately, but with new probabilities $\tilde{\alpha}_{c_i}$, since they are thinned independently by impatience.

\begin{figure}[htb]
   \begin{center}
    \includegraphics[scale=0.40]{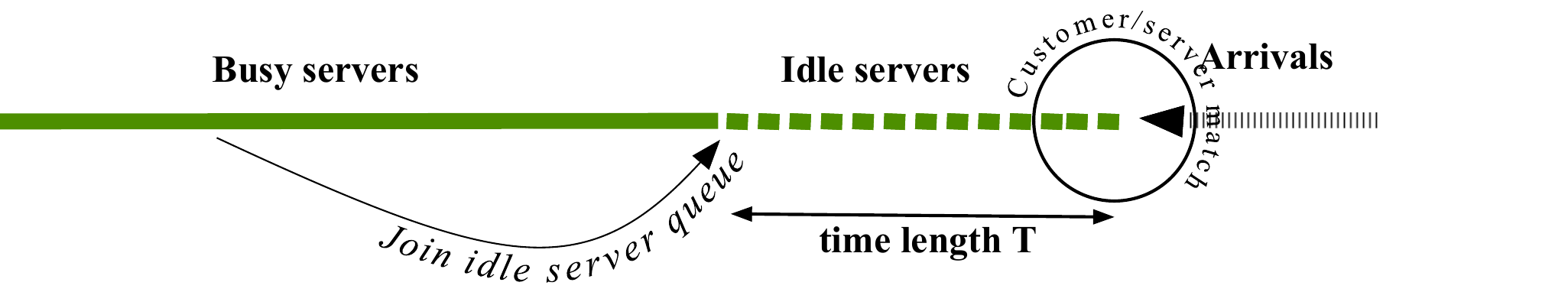}
   \end{center}
   \caption{Operation of system in QD mode}\label{fg2}
\end{figure}

\begin{figure}[htb]
   \begin{center}
    \includegraphics[scale=0.40]{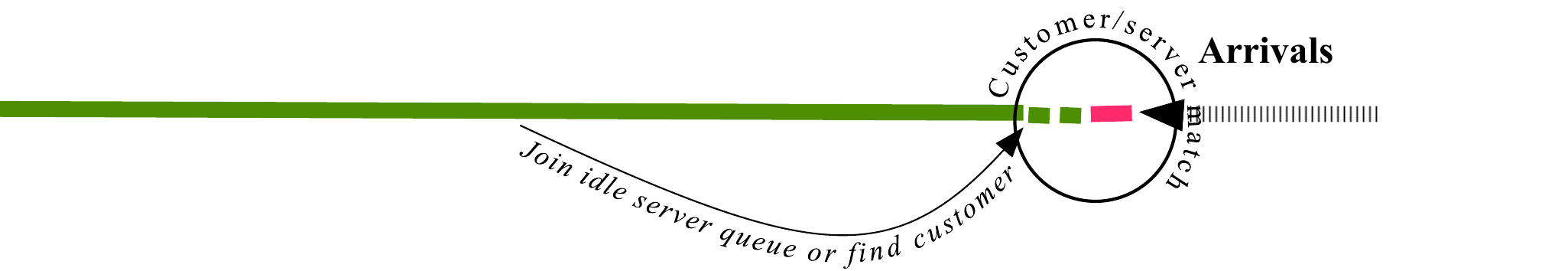}
   \end{center}
\caption{Operation of system in QED mode}\label{fg3}
\end{figure}

In QD mode, there is a queue of idle servers, each server, on completing a service, joins the end of this queue (see Figure \ref{fg3}).  A server reaches the head of the queue after an idle time $T$, and matches with the first compatible customer.
Customers never wait and are of i.i.d. types with probabilities $\alpha_{c_i}$.

In QED mode, the system alternates infrequently between periods with a short queue of waiting customers and periods with a short queue of idle servers (see Figure 3).  All servers are almost always busy, customers immediately enter service or wait a short time, and abandonments are rare.

The key assumption necessary for the matching rates to be according to the FCFS infinite bipartite matching model is that approximately the sequence of customers entering service has i.i.d. types and the sequence of servers that become available and that start service has i.i.d. types.

When there is no resource pooling, the system decomposes into subsystems as stated in (\ref{eqn.decomposition}).   To illustrate this situation, assume that the staffing is such that the system decomposes into 3 sub-systems.  Sub-system $(\C^{(1)},\S^{(1)})$ has the least staffing, and so will receive the worst service, while sub-system $(\C^{(3)},\S^{(3)})$ has the most abundant staffing, and will receive the best service, with subsystem $(\C^{(2)},\S^{(2)})$ inbetween.

With no abandonments, starting with small $\lambda$, the system will be stable under FCFS-ALIS.  If we let $\lambda$  increase, servers in $\S^{(1)}$ will become fully utilized, and the queue of customers of types in  $\C^{(1)}$ will become unstable, while $(\C^{(2)},\S^{(2)}), (\C^{(3)},\S^{(3)})$ remain stable.  As $\lambda$  increases further, servers $\S^{(2)}$ will also become fully utilized, and the queue of customers of types $\C^{(2)}$ will also become unstable, while $(\C^{(3)},\S^{(3)})$ remain stable, and finally when $\lambda$ becomes even larger, all the servers will be fully utilized and
the queues of all types will become unstable.

Furthermore, once the whole system is unstable, the longest waiting customers will all be of types in $\C^{(1)}$, waiting for servers $\S^{(1)}$.
Servers $\S^{(2)}$ will skip those longest waiting customers and will serve customers of types $\C^{(2)}$ which will have shorter waits, and servers $\S^{(3)}$ will skip all the waiting customers of type $\C^{(1)},\C^{(2)}$, and serve customers of types $\C^{(3)}$, which will have the shortest waits.

Under abandonments, for high enough $\lambda$, customers of types $\C^{(1)}$ will have the highest fraction of abandonments and the longest wait, customers of types $\C^{(3)}$ will have the smallest fraction of abandonments and the shortest wait, and types  $\C^{(2)}$ will be inbetween.

On the other hand, if the system is stable, for small enough $\lambda$, under ALIS servers of types $\S^{(1)}$ will have shorter idle periods than servers $\S^{(2)}$, who in turn will have shorter idle periods than servers $\S^{(3)}$.

\begin{figure}[htb]
 \begin{center}
    \includegraphics[scale=0.40]{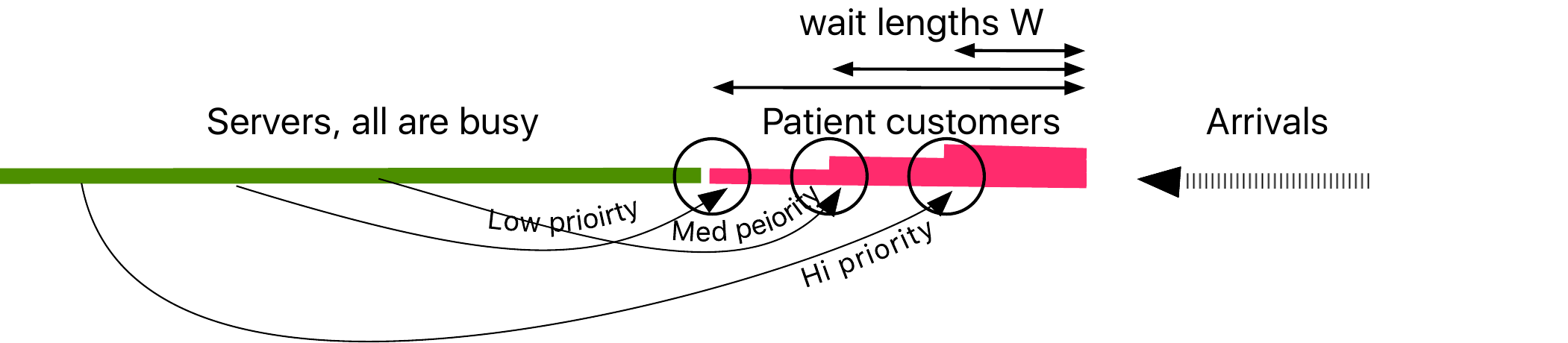}
   \end{center}
 \begin{center}
    \includegraphics[scale=0.40]{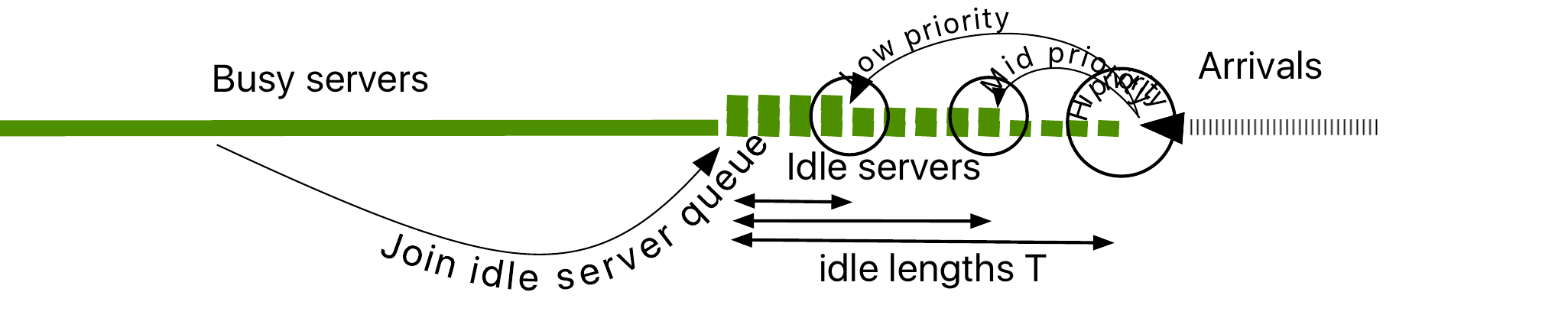}
   \end{center}
 \begin{center}
    \includegraphics[scale=0.40]{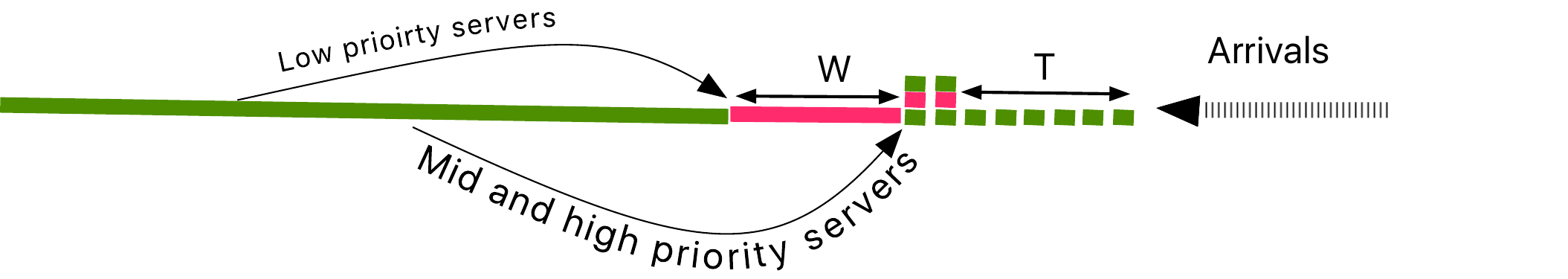}
   \end{center}
   \caption{Operation of system in ED mode (top), QD mode (middel) and QED mode (bottom)}\label{fg4}
\end{figure}

Figure \ref{fg4} shows how such a system will behave under our conjecture.

The three illustrations depict behavior under ED (top figure), QD (middle figure), and in the bottom one, $(\C^{(1)},\S^{(1)})$ is in ED, $(\C^{(2)},\S^{(2)})$ is in QED and $(\C^{(3)},\S^{(3)})$ is in QD.   In ED mode, customers of types $\C^{(1)}$ wait the longest time, before being served by servers of types $\S^{(1)}$.  Servers of types $\S^{(2)}$ skip customers of type $\C^{(1)}$ and serve customers of type $\C^{(2)}$ that have a shorter wait, and servers of types $\S^{(3)}$ skip customers of type $\C^{(1)},\C^{(2)}$ and serve customers of type $\C^{(3)}$ that have the shortest wait.
In  QD mode, servers of the high priority customers have longer idle periods, and are available for the high priority customers immediately, while servers of lower priority customers have shorter idle times, and are further in the queue of idle servers, and their customers will skip the high priority servers in the queue, which are incompatible with them.   In the mixed mode, low priority customers will wait, while high priority customers will have a queue of idle servers ready to serve them.
The midlevel customers and their servers will have  periods of short waits alternating with periods of short idle times.


\section{Design Algorithms}
\label{sec.algorithm}

\subsection{General Strategy}
The setup for the design problem is given by server and customer types and the compatibility graph.  To these  are added the patience distributions of the customers, and the service time distributions for each customer/server type pair.  Of the latter only the average service time and service rate are required for the design.  Next the data includes the arrival rate $\lambda$, and the customer type frequencies $\alpha_{c_i}$.

The first design decision is whether we wish to provide uniform service levels to all types, or whether we divide them to classes of varying priorities.  In the first case we design the system to have complete resource pooling, in the latter case we partition the customers, and assign subsets of servers to each priority class in such a way that the different priorities will receive different service levels under FCFS-ALIS (details follow later).

Next we need to decide on mode of operation, ED, QD or QED, and on parameters of quality of service and of utilization.  Consider the case of complete resource pooling. If we decide on ED, then quality of service will be determined by specifying the  average waiting time uniformly for all types of customers (which also determines fractions of abandonments).  If we decide on QD, then level of utilization will be determined by specifying an average idle time uniform for all types of servers.  If we decide on QED the required system will need to have servers almost fully utilized and zero or short customer queues.  In the case of decomposition to several classes with different priorities, mode of operation and quality parameters need to be determined for each sub-class separately.

For given $\lambda$ and $\alpha_{c_i}$, having specified the degree of decomposition and the quality parameters, there will be many staffing combinations of servers that will satisfy these requirements.  The next design decision will specify which of these staffing designs we are to choose.

As we stated in the introduction, we can do this in two different ways:  We can pre-specify the fraction of the total number of services $\beta_{s_j}$ which each type of server performs, or we can pre-specify the fraction $\theta_{s_j}$ of the total number of servers for each type.  The first is easier, and we will describe the algorithms based on pre-specified $\beta_{s_j}$, in which we calculate the required $n$ and $n_{s_j}$, and determine $\theta_{s_j}$.  If $\theta_{s_j}$ are pre-specified we will need to solve numerically for the $\beta_{s_j}$ that will yield these
values of $\theta_{s_j}$.   We outline the numerical procedure as well.

Once we determined decomposition level and quality parameters, modified $\tilde{\alpha}_{c_j}$ can be determined and with pre-specified $\beta_{s_j}$ we  then   use the bipartite infinite matching model, formula (\ref{eqn.matches}), to obtain the matching rates $r_{c_i,s_j}$.
Once we have the matching rates, we can calculate the amount of work required from each type of server, and this determines, by Little's law, the number of servers that are needed of each type in order to meet the requested quality of service and utilization.

In the following sub-sections we show how to perform these steps for complete resource pooling under each of the three regimes, as well as for the case of systems decomposed to several sub-systems with differing priorities.  We then outline the numerical calculations for the case that $\theta_{s_j}$ is pre-specified in Sub-section \ref{sec.newton}.

We illustrate the calculations for several examples in Section \ref{sec.examples}, and demonstrate the effectiveness of the heuristics through simulation.

\subsection{Design for Quality Driven Service}
Here the traffic intensity is $<1$, and customers almost never wait, and therefore even more rarely abandon.  There are almost always some idle servers waiting for customers, and because of ALIS,   servers of different types all have  the same idle time distribution.   The quality parameter in this case is the value $T$ of the average idle time.  It is a measure of the utilization of the servers.
Because there are virtually no abandonments, the patience time distribution is not required as input.

\begin{center}
\begin{tabular}{|p{0.92\textwidth}|}
\hline
\hfill\textbf{Algorithm for QD}\hfill\mbox{}\\
\hline
\noindent{\bf Input:}
\begin{compactitem}
\item Compatibility graph $\mathcal{G}$
\item Arrival rate $\lambda$
\item Fractions of customer types $\alpha_{c_i}$
\item Mean service times $m_{c_i,s_j}$
\end{compactitem}\\
\hline
\noindent{\bf Requested quality of service parameter:}
\begin{compactitem}
\item Mean server idle time after each service $T$
\end{compactitem}\\
\hline
\noindent\textbf{Design parameters:}
\begin{compactitem}
\item Fraction of services performed by each server type $\beta_{s_j}$
\end{compactitem}\\
\hline
\noindent\textbf{Algorithm}:\\
\parbox{0.8\linewidth}{
\begin{align*}
&\textrm{\bf Check $\alpha,\beta$ for complete resource pooling}\\
&\textrm{Compute matching rates} &&& r_{c_i,s_j} &:= \textrm{use Equation }\eqref{eqn.matches} \\
&\textrm{Compute staffing levels}&&& n_{s_j} &:=  \sum_{c_i \in \C(s_j)}  \lambda r_{c_i,s_j} \big( m_{c_i,s_j} + T\big)
\end{align*}}\\
\hline
\noindent{\bf Output:}\\
\begin{compactitem}
\item Required workforce $n_{s_j}$
\end{compactitem}\\
\hline
\end{tabular}
\end{center}

\subsection{Design for Efficiency Driven Service}
Here the traffic intensity is $>1$, servers are always busy and customers always need to wait, and a certain fraction will abandon.  By FCFS,  customers of different types all have  the same waiting  time distribution, and the system demonstrates global FCFS (this term was coined by Talreja and Whitt \cite{talreja-whitt:07}).
The system is stabilized by abandonments, with average waiting time $W$.   This means that approximately a fraction $F_{c_i}(W)$ of customers of type $c_i$ will abandon.    The value of $W$ is the quality of service parameter here, so that customers with patience less than $W$ do not get served, while customers with patience that exceeds $W$ get served after a wait of $W$.
Since customers are thinned independently by impatience, we need to calculate the total effective arrival rate (of patient customers) and we need to adjust the fractions $\alpha_{c_i}$ of customer types entering service.

\begin{center}
\begin{tabular}{|p{0.92\textwidth}|}
\hline
\hfill\textbf{Algorithm for ED}\hfill\mbox{}\\
\hline
\noindent{\bf Input:}
\begin{compactitem}
\item Compatibility graph $\mathcal{G}$
\item Arrival rate $\lambda$
\item Fractions of customer types $\alpha_{c_i}$
\item Patience distributions $F_{c_i}(\cdot)$
\item Mean service times $m_{c_i,s_j}$
\end{compactitem}
\\\hline
\noindent{\bf Requested quality of service parameter:}
\begin{compactitem}
\item Mean waiting time $W$
\end{compactitem}
\\\hline
\noindent\textbf{Design parameters:}
\begin{compactitem}
\item Fraction of services performed by each server type $\beta_{s_j}$
\end{compactitem}
\\\hline
\noindent\textbf{Algorithm:}\\
\parbox{\linewidth}{
\begin{align*}
&\textrm{Compute expected fraction of abandonments} &&& p_{c_i} &:= F_{c_i} (W)\\
&\textrm{Compute the total effective arrival rate} &&& \tilde{\lambda}&:=\sum_{i=1}^I \alpha_{c_i}\lambda (1-p_{c_i})\\
&\textrm{Adjust fraction of each customer type}&&& \tilde{\alpha}_{c_i} &:= \frac{\alpha_{c_i}\lambda  (1-p_{c_i})}{\tilde{\lambda}}\\
&\textrm{\bf Check $\tilde{\alpha},\beta$ for complete resource pooling}\\
&\textrm{Compute matching rates} &&& r_{c_i,s_j} &:= \textrm{use Equation }\ \eqref{eqn.matches} \\
&\textrm{Compute staffing levels}&&& n_{s_j} &:=  \sum_{c_i \in \C(s_j)}  \tilde{\lambda} r_{c_i,s_j}
m_{c_i,s_j}
\end{align*}
}
\\\hline
\noindent{\bf Output:}\\
\begin{compactitem}
\item Required workforce $n_{s_j}$
\end{compactitem}
\\\hline
\end{tabular}
\end{center}

\subsection{Design for Quality and Efficiency Driven Service}
Given the arrival rates, there is a unique FCFS system that will supply QED service, with servers almost always busy, most customers either don't wait or wait a short time, and few abandonments.
The calculation of QED design follows the same steps as for QD with $T=0$ and for ED with $W=0$.

%

\subsection{Design for Differentiated Service}

We now consider the case where we would like to give customers graded service levels, from high priority to  standard priority to low priority customer types.
We have a partition of customer types $\C$  into $\C^{(1)},\ldots,\C^{(L)}$, where customers of types $c_i\in \C^{(l)}$ have higher priority than customers of type  in $\C^{(l-1)}$, $l=2,\ldots,L$.
Customers in class $\C^{(l)}$ will then have a set of servers $\S^{(l)}$, so that each customer type $c_i\in \C^{(l)}$ will have at least one compatible server type $s_j \in \S^{(l)}$.  In this decomposition
to subsystems $(\C^{(l)},\S^{(l)})$, we will allow $s_j \in \S^{(l)}$ to serve $c_i \in \C^{(k)},\,k \ge l$, but will not allow $s_j \in \S^{(l)}$ to serve $c_i \in \C^{(k)},\,k<l$. In other words, we redesign the compatibility graph $\mathcal{G}$ by eliminating all links from $\S^{(l)}$ to  $\C^{(k)}$ for $k<l$, but we preserve the
links to higher priority customers in $\C^{(k)}$ for $k>l$, since when we use FCFS, these links will hardly ever be used, because servers in $\S^{(l)}$ will be behind all servers in $\S^{(k)}$ for $k>l$ almost all the time.

The priorities will be translated into quality of service parameters:  classes $\C^{(1)},\ldots,\C^{(l)}$ will be served in ED mode with $W_1>\cdots>W_l\ge 0$, classes $l+1,\ldots,L$ will be served in QD mode, with $0< T_{l+1}<\cdots<T_L$.  Class $l$ may be in QED mode.

\begin{center}
\begin{tabular}{|p{0.92\textwidth}|}
\hline
\hfill\textbf{Algorithm for Differentiated Service}\hfill\mbox{}\\
\hline
\noindent{\bf Input:}
\begin{compactitem}
\item Compatibility graph $\mathcal{G}$
\item Arrival rate $\lambda$
\item Fractions of customer types $\alpha_{c_i}$
\item Patience distributions $F_{c_i}(\cdot)$
\item Mean service times $m_{c_i,s_j}$
\end{compactitem}\\
\hline
\noindent{\bf Requested quality of service parameters:}
\begin{compactitem}
\item Partition of customer types by priority into $\C^{(1)},\ldots,\C^{(L)}$
\item Quality of service parameters:  $W_1 > \cdots > W_l=0=T_l < T_{l+1} < \cdots < T_L$
\end{compactitem}
\\\hline
\noindent\textbf{Design parameters:}
\begin{compactitem}
\item Choose partition of server types $\S^{(1)},\ldots,\S^{(L)}$
\item Eliminate links from $\S^{(l)}$ to $\C^{(k)}$, for $k<l$
\item Assign fraction of services performed by each server type $\beta_{s_j}$, within $\S^{(l)}$
\end{compactitem}\\
\hline
\noindent\textbf{Algorithm}:\\
\begin{compactitem}
\item For subsystem $(\C^{(l)},\S^{(l)})$, $l=1,\ldots,L$:
\item Apply appropriate design algorithm for ED, QED, or QD to the subsystem
\end{compactitem}\\
\hline
\noindent{\bf Output:}\\
\begin{compactitem}
\item Redesigned compatibility graph $\mathcal{G}$
\item Required workforce $n_{s_j}$
\end{compactitem}\\
\hline
\end{tabular}
\end{center}

\subsection{Numerical calculations when fraction of servers of each type is specified}
\label{sec.newton}
In the previous sub-sections we illustrated how to obtain the vector of values $(n_{s_1},\ldots,n_{s_J})$ for given values of $(\beta_{s_1},\ldots,\beta_{s_J})$, which also determines $(\theta_{s_1},\ldots,\theta_{s_J})$.  We therefore have a function $H:\mathbb{R}^J \to \mathbb{R}^J$, which calculates $(\theta_{s_1},\ldots,\theta_{s_J}) = H(\beta_{s_1},\ldots,\beta_{s_J})$.  We let $H^j(\beta_{s_1},\ldots,\beta_{s_J})$ denote the $j$th element of $H$.   We now need to perform the inverse calculation, of $(\beta_{s_1},\ldots,\beta_{s_J}) = H^{-1}(\theta_{s_1},\ldots,\theta_{s_J})$.  We do not know whether $H$ is one-to-one, so $H^{-1}$ may be multi-valued, in which case we would like to find just one inverse vector.   To obtain such a vector we consider, for given $\theta_{s_1},\ldots,\theta_{s_J}$ and some proposed $\beta_{s_1},\ldots,\beta_{s_J}$, the squared sum of differences:
\begin{equation}
\Delta =   \sum_{j=1,\ldots,J} \big(H^j(\beta_{s_1},\ldots,\beta_{s_J}) - \theta_j\big)^2 .
\label{eqn:Delta}
\end{equation}
For given $(\theta_{s_1},\ldots,\theta_{s_J})$ we wish to find $(\beta_{s_1},\ldots,\beta_{s_J})$ that will minimize $\Delta$ or that will solve $\Delta=0$.

Numerical solution of this problem can be performed by a wide choice of  minimization software.  We illustrate the solution for the examples of Sections \ref{sec.example1} and \ref{sec.example3}, in Section \ref{sec.minimization}.


\section{Examples of Designs and Simulation Results}
\label{sec.examples}

In this section we describe three examples, for  each of which we have prepared several designs, under several modes of operation, and  assuming Poisson arrivals using a range of values for $\lambda$. We have restricted the examples to Poisson arrivals to keep the presentation manageable and also because experiments confirm that the same conclusions remain valid for general renewal streams. We have then performed extensive simulation runs on each of these designs.
Our purpose in this section is threefold:
\begin{itemize}
\item[(i)]
Illustrate the implementation of the algorithms;
\item[(ii)]
Examine the validity of the matching rates conjecture;
\item[(iii)]
Evaluate the efficacy of our designs.
\end{itemize}

The first example system has 3 customer types and 3 server types with an almost complete bipartite compatibility graph.  We have designed operation of this system with complete resource pooling, in ED, QD and QED mode.  The purpose of this example is to assess pooled service designs.  This example is important, because it is for this system topology that Foss and Chernova \cite{foss-chernova:98} have shown that calculation of exact matching rates depends of the full shape of processing time distributions and is intractable.

The second example has 5 customer types and 5 server types with a Hamiltonian bipartite compatibility graph.  This example is used to assess differentiated service, with the customer types divided into high, standard and low priority customers.
Matching rates for the decomposed system are easy to obtain.  What we show however is that the system achieves the desired graded service specifications under FCFS-ALIS.

The third example has 6 customer types and 6 server types with a symmetric compatibility graph that has degree 3 for all nodes.  The purpose of considering this example is to examine validity of the matching conjecture in a  complex graph.

The  designs depend on the service rates for each link, but not on the actual distributions of service times.  To examine  the validity of the matching conjecture and to assess the efficacy of the designs, we have chosen to simulate service time distributions which are very different, including uniform in a finite range, exponential and Pareto.

Our main conclusions from the simulations of these examples are:
\begin{itemize}
\item[(i)]
The matching rates conjecture seems to be valid under ED, QD and QED mode, for the
whole range of $\lambda$ values, and under all the different distributions of service times.  For small values of $\lambda$ the deviations are slightly larger. This can be partly explained by the fact that the algorithm yields real numbers for  $n_{s_j}$, while in the implemented design rounded integer values are used.
\item[(ii)]
In ED mode, for large values of $\lambda$ we get convergence to the exact values of $W$  with very small variability in waiting times, and exact abandonment rates.  Similarly under QD, for large values of $\lambda$ we get convergence to the exact values of $T$  with very small variability in idle times, and almost all customers are not waiting for service.
\item[(iii)]
Most important, it seems that for small values of $\lambda$, while waiting times in ED mode and idle times in QD mode are quite variable, the average waiting time in ED and the average idle time in QD are almost exactly as designed.  This indicates that our design heuristic may be effective already for a moderate number of servers.
\item[(iv)]
Convergence in the QED mode is not appreciably worse than in the ED or QD modes.
\item[(v)]
The system with  differentiated service performs as designed, under FCFS-ALIS.
\item[(vi)] The results do not seem to depend on the service time distributions.
\end{itemize}

We now present the three examples with detailed simulation results.  The reported simulation results for each design have been obtained as the average of 1,000 runs, where each run consists of 1,250,000 customers. However, the first 250,000 customers have been removed from the results to account for a possible startup effect.

In Sections \ref{sec.example1}-\ref{sec.example3} we first calculate designs with specified service fractions $\beta_{s_j}$, and finally, in Section \ref{sec.minimization} we demonstrate calculation of designs with specified $\theta_{s_j}=n_{s_j}/n$ rather than specified $\beta_{s_j}$.

\subsection{ Example 1 -- \boldmath $3\times3$ Almost Complete Graph with Pooled Service}
\label{sec.example1}

In this example we investigate pooled service designs. The system is specified below, where Exp$(a)$ denotes the exponential distribution with rate $a$, U$(a,b)$ is the uniform distribution on the interval $(a,b)$ and Pareto$(k,a)$ is the Pareto distribution $F(t) = 1 - (k/t)^a$ for $t > k$.

\begin{center}
\begin{tabular}{|p{0.92\textwidth}|}
\hline
\hfill\textbf{ Example 1 -- System and Data}\hfill\mbox{}\\
\hline
There are 3 types of customers and 3 types of servers. The total arrival rate is parameterized by $\lambda$,  The graph and the values of $\alpha_{c_i} \lambda$ are described in the following figure:
   \begin{center}
    \includegraphics[scale=0.30]{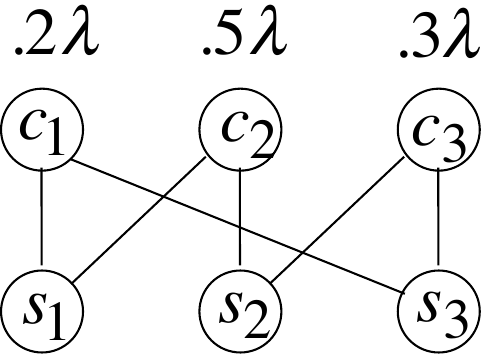}
   \end{center}
The patience times and service time distribution are given in the tables below.
\begin{center}
\begin{tabular}{|c|c|}
\hline
\multicolumn{2}{|c|}{Patience time distributions}\\
\hline
& $F_{c_i}$    \\
\hline
$c_1$ & Exp(0.1)\\
$c_2$ & U(0,10) \\
$c_3$ & Exp(0.2)\\
\hline
\end{tabular}
\\
\begin{tabular}{|c|ccc|}
\hline
\multicolumn{4}{|c|}{Service time distributions}\\
\hline
$G_{c_i,s_j}$    & $c_1$ & $c_2$ & $c_3$ \\
\hline
$s_1$ & Pareto(2, 3) & Exp(0.125)  &\\
$s_2$ & & Exp(0.2) & U(2, 6) \\
$s_3$ & Pareto(3, 3) & & U(1, 5) \\
\hline
\end{tabular}
\quad \begin{tabular}{|c|ccc|}
\hline
\multicolumn{4}{|c|}{Mean service times}\\
\hline
$m_{c_i,s_j}$    & $c_1$ & $c_2$ & $c_3$ \\
\hline
$s_1$ & 3 & 8  &\\
$s_2$ & & 5 & 4 \\
$s_3$ & 4.5 & & 3 \\
\hline
\end{tabular}

\end{center}
Only the \emph{mean} service times are used by the design algorithms. The full distributions are used in the simulations.\\
\hline
\end{tabular}
\end{center}

In the designs for Example 1 we take as service fractions:
$\beta_{s_1}=0.3$, $\beta_{s_2}=0.3$, $\beta_{s_3}=0.4$.\\

\noindent {\bf ED design:}  We specify the average waiting time $W=1$, corresponding to approximately 25\% of the mean service times.
For the given patience distributions this entails abandonment rates of approximately 10\% for  customers of types 1 and 2, and of 18\% for customers of type 3.
We calculate the effective arrival rates of customers that do get served after a wait of $W=1$:
\[
1-F_{c_1}(W) = e^{-0.1 W} = 0.905,  \quad 1-F_{c_2}(W) = (10-W)/10=0.9,
\quad 1-F_{c_3}(W) = e^{-0.2 W} = 0.819,
\]
so
\[
\lambda_{c_1}(1-F_{c_1}(W)) = 0.2 \lambda \times 0.905 = 0.181\lambda,\quad
\lambda_{c_2}(1-F_{c_2}(W)) = 0.450 \lambda, \quad \lambda_{c_3}(1-F_{c_3}(W))=0.246\lambda,
\]
and thus the effective arrival rate equals
\[
\tilde{\lambda} = (0.181+0.450+0.246) \lambda = 0.877 \lambda .
\]
Hence, the adjusted values of $\alpha_{c_j}$ are:
\[
\tilde{\alpha}_{c_1} = \frac{0.18}{0.88}= 0.206,\quad
\tilde{\alpha}_{c_2} = \frac{0.45}{0.88}= 0.513,\quad
\tilde{\alpha}_{c_3} = \frac{0.25}{0.88}= 0.281.
\]

\noindent {\bf QD design:}  We take an average idle time of $T=0.5$.  This corresponds to an utilization of approximately 0.9.

\noindent {\bf QED design:}  The unadjusted values of $\lambda, \alpha_{c_i},\,m_{c_i,s_j}$ are used.\\

It is readily verified that in all three regimes (ED, QD and QED), Conditions (\ref{eqn.crp}) are satisfied, so complete resource polling holds.
From the algorithms we obtain the {\em calculated required workforce} for the three designs shown in Table \ref{tbl:rw1}.

\begin{table}[htb]
\begin{center}
\begin{tabular}{|c|ccc|ccc|ccc|}
\hline
 & \multicolumn{9}{|c|}{Required workforce}\\
 \hline
& \multicolumn{3}{|c|}{ED regime} & \multicolumn{3}{|c|}{QED regime} & \multicolumn{3}{|c|}{QD regime} \\
\hline
$\lambda$
  & $n_{s_1}$ & $n_{s_2}$ & $n_{s_3}$
  & $n_{s_1}$ & $n_{s_2}$ & $n_{s_3}$
  & $n_{s_1}$ & $n_{s_2}$ & $n_{s_3}$
  \\
\hline
20   &   39 & 25 & 25     &  44 & 29 & 29      &    47 & 32 & 33   \\
40   &   77 & 51 & 51     &  88 & 58 & 57      &    94 & 64 & 65   \\
60   &  116 & 76 & 76     &  131 & 87 & 86     &    140 & 96 & 98   \\
100  &  194 & 127 & 127   &  219 & 144 & 144   &    234 & 159 & 164 \\
200  &  387 & 254 & 255   &  438 & 288 & 287   &    468 & 318 & 327 \\
\hline
\end{tabular}
\caption{Calculated required workforce for Example 1}\label{tbl:rw1}
\end{center}
\end{table}

The simulation results for Example 1 are listed in the tables below. We note that the histograms below only depict the waiting times and idle times greater than zero (so probability mass at zero is not shown).

\begin{center}
\begin{tabular}{|c|ccc|ccc|ccc|}
\hline
 & \multicolumn{9}{|c|}{Matching rates}\\
\hline
& \multicolumn{3}{|c|}{ED regime} & \multicolumn{3}{|c|}{QED regime} & \multicolumn{3}{|c|}{QD regime} \\
\hline
 & \multicolumn{9}{|c|}{Theoretical}\\
\hline
$r_{c_i,s_j}$    & $c_1$ & $c_2$ & $c_3$  & $c_1$ & $c_2$ & $c_3$& $c_1$ & $c_2$ & $c_3$\\
\hline
$s_1$ &  0.038 & 0.262 &       &  0.042 & 0.258 &       &  0.042 & 0.258 &         \\
$s_2$ &        & 0.251 & 0.049 &        & 0.242 & 0.058 &        & 0.242 & 0.058   \\
$s_3$ &  0.168 &       & 0.232 &  0.158 &       & 0.242 &  0.158 &       & 0.242   \\
\hline
 & \multicolumn{9}{|c|}{$\lambda=20$}\\
\hline
$r_{c_i,s_j}$    & $c_1$ & $c_2$ & $c_3$  & $c_1$ & $c_2$ & $c_3$& $c_1$ & $c_2$ & $c_3$\\
\hline
$s_1$ &   0.046 & 0.262 &       &    0.048 & 0.260 &       &  0.047 & 0.258 &       \\
$s_2$ &         & 0.241 & 0.056 &          & 0.239 & 0.064 &        & 0.241 & 0.065 \\
$s_3$ &   0.164 &       & 0.230 &    0.155 &       & 0.234 &  0.153 &       & 0.236 \\
\hline
 & \multicolumn{9}{|c|}{$\lambda=60$}\\
\hline
$r_{c_i,s_j}$    & $c_1$ & $c_2$ & $c_3$  & $c_1$ & $c_2$ & $c_3$& $c_1$ & $c_2$ & $c_3$\\
\hline
$s_1$ &   0.041 & 0.261 &       &   0.045 & 0.258 &       & 0.045 & 0.257 &       \\
$s_2$ &         & 0.248 & 0.051 &         & 0.242 & 0.061 &       & 0.243 & 0.062 \\
$s_3$ &   0.167 &       & 0.232 &   0.156 &       & 0.237 & 0.155 &       & 0.238 \\
\hline
 & \multicolumn{9}{|c|}{$\lambda=200$}\\
\hline
$r_{c_i,s_j}$    & $c_1$ & $c_2$ & $c_3$  & $c_1$ & $c_2$ & $c_3$& $c_1$ & $c_2$ & $c_3$\\
\hline
$s_1$ &   0.039 & 0.261 &       &  0.043 & 0.259 &       &  0.043 & 0.258 &       \\
$s_2$ &         & 0.250 & 0.049 &        & 0.242 & 0.059 &        & 0.242 & 0.059 \\
$s_3$ &   0.168 &       & 0.232 &  0.157 &       & 0.240 &  0.157 &       & 0.241 \\
\hline
\end{tabular}
\end{center}

\begin{center}
\begin{tabular}{|c|c|c|c|}
\hline
\multicolumn{4}{|c|}{Customer waiting times and server idle times}\\
\hline
ED Regime & QED Regime & QED Regime & QD Regime\\
Waiting times & Waiting times & Idle times & Idle times\\
\hline
&&&\\[-2mm]
\includegraphics[width=0.24\textwidth]{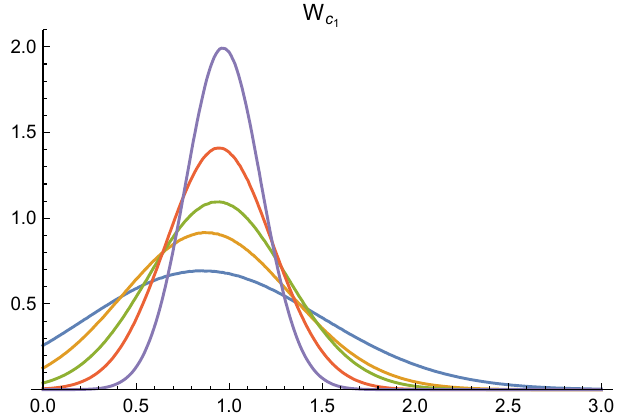} & \includegraphics[width=0.24\textwidth]{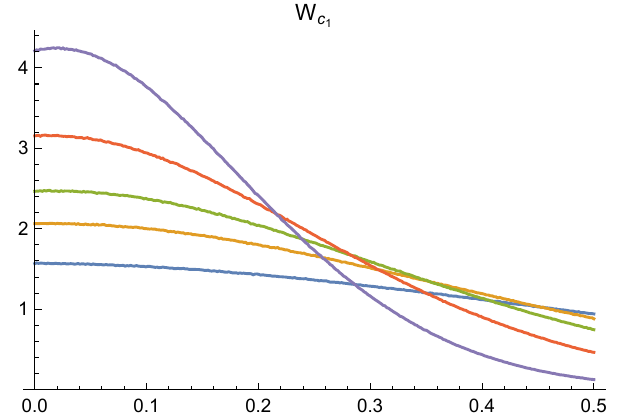} &\includegraphics[width=0.24\textwidth]{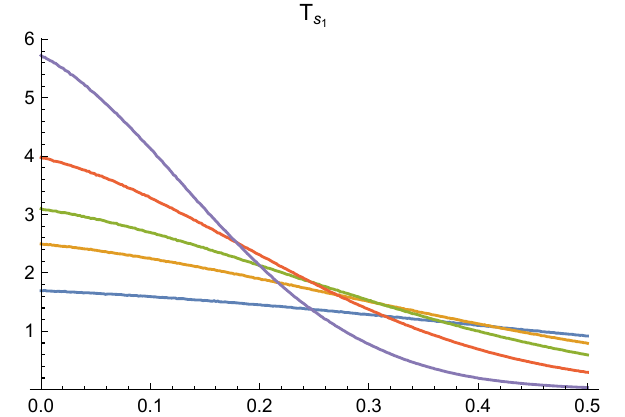} &\includegraphics[width=0.24\textwidth]{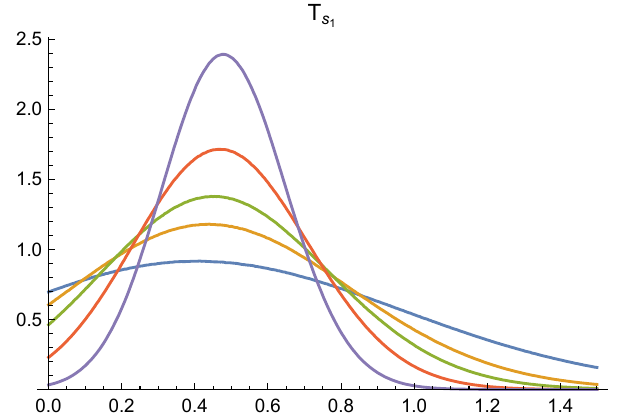} \\
\includegraphics[width=0.24\textwidth]{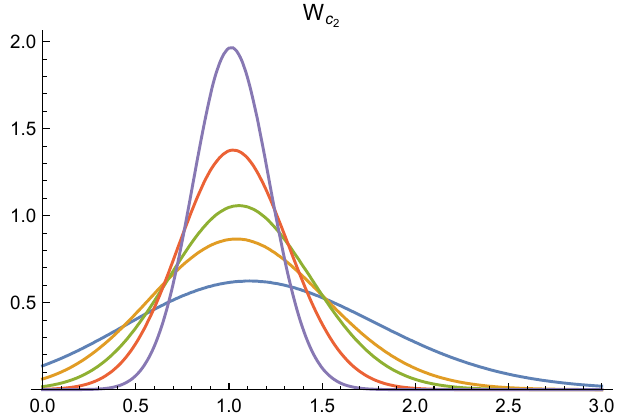} & \includegraphics[width=0.24\textwidth]{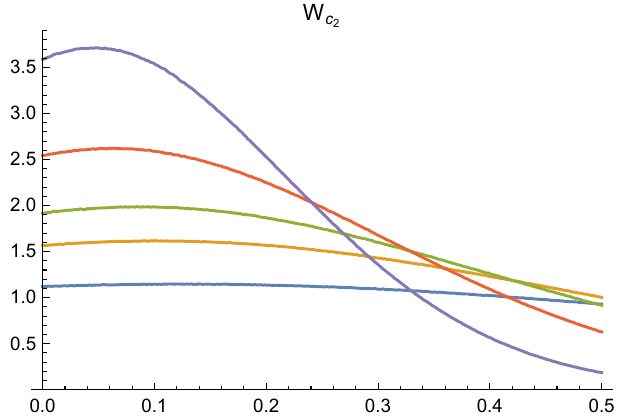} &\includegraphics[width=0.24\textwidth]{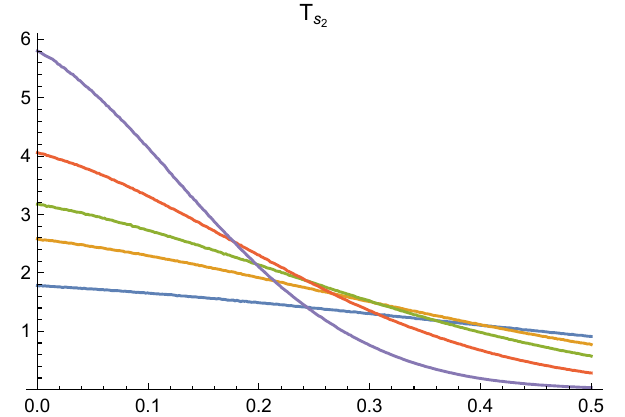} &\includegraphics[width=0.24\textwidth]{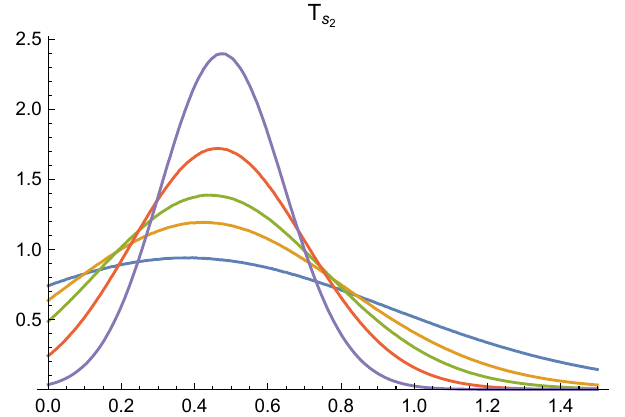}   \\
\includegraphics[width=0.24\textwidth]{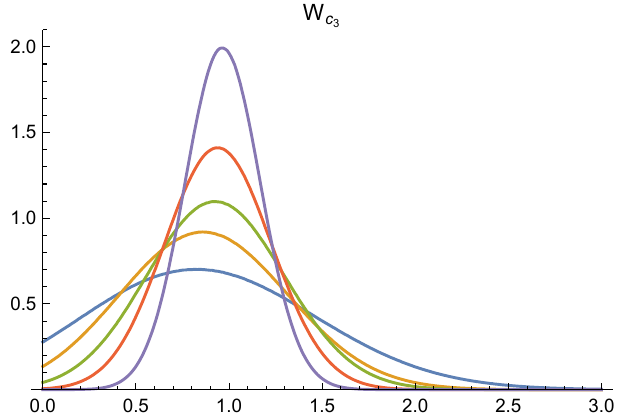} & \includegraphics[width=0.24\textwidth]{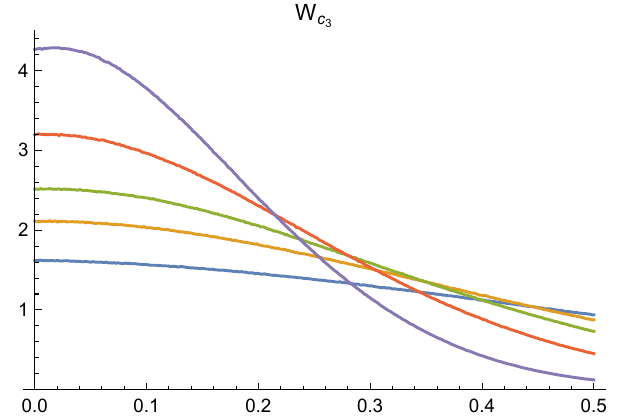}  &\includegraphics[width=0.24\textwidth]{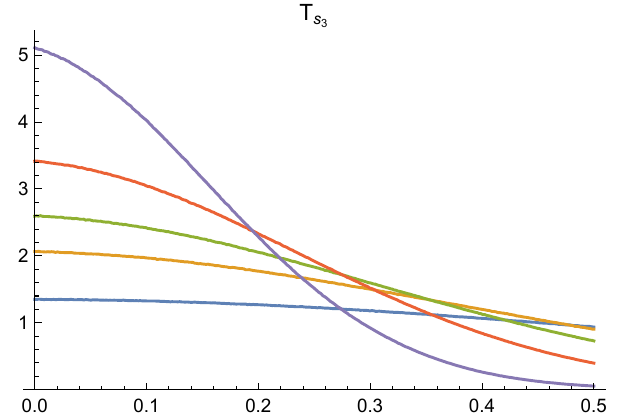} &\includegraphics[width=0.24\textwidth]{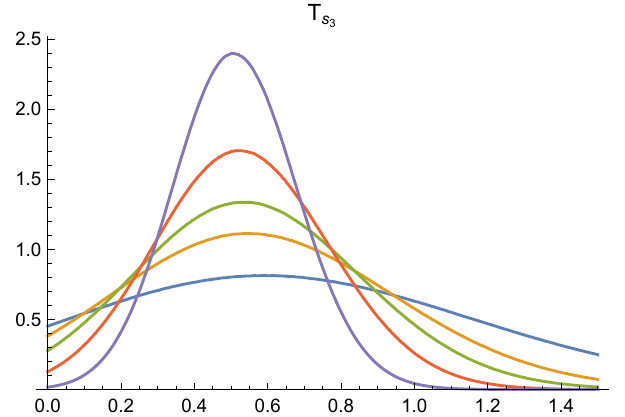}     \\
\hline
\multicolumn{4}{|c|}{\includegraphics[width=0.6\textwidth]{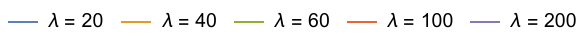}}\\
\hline
\end{tabular}
\end{center}


\begin{center}
\begin{tabular}{|l|c|c|c|c|c|c|}
\hline
& \multicolumn{6}{|c|}{Fraction of no wait and of no idling}\\
\hline
    & \multicolumn{2}{|c|}{ED regime}&    \multicolumn{2}{|c|}{QED regime}& \multicolumn{2}{|c|}{QD regime} \\
\hline
$\lambda$    & No waiting & No idling & No waiting &  No idling  & No waiting   & No idling \\
\hline
20   & 0.047 & 0.946       & 0.444 &  0.540      & 0.814    &  0.181         \\
60   & 0.003 & 0.997       & 0.420 &  0.571      & 0.947    &  0.053           \\
200  & 0.000 & 1.000       & 0.410 &  0.585      & 0.999    &  0.001             \\
\hline
\end{tabular}
\end{center}

\begin{table}[h]
 \begin{center}
\begin{tabular}{|c|ccc|ccc|ccc|}
 \hline
&  \multicolumn{9}{|c|}{Abandonment rates}  \\
\hline
&  \multicolumn{3}{|c|}{ED regime}  &  \multicolumn{3}{|c|}{QED regime}&  \multicolumn{3}{|c|}{QD regime}\\
\hline
$\lambda$
  & $c_1$ & $c_2$ & $c_3$
  & $c_1$ & $c_2$ & $c_3$
  & $c_1$ & $c_2$ & $c_3$
  \\
\hline
20 &  0.089 & 	0.123	 & 0.168  & 0.020 & 	0.035 & 	0.039 & 0.003 & 	0.008	 & 0.006\\
60 &  0.091 & 	0.109	 & 0.173  & 0.014 & 	0.020 & 	0.027 & 0.000 & 	0.001	 & 0.001\\
200 & 0.093 & 	0.102	 & 0.177  & 0.008 & 	0.010 & 	0.016 & 0.000 & 	0.000	 & 0.000\\
\hline
Design & 0.095	&0.100	&0.181 & 0.000 & 0.000 & 0.000 & 0.000 & 0.000 & 0.000 \\
\hline
\end{tabular}
\end{center}
\caption{Simulation results for Example 1}\label{tbl:res1}
\end{table}

The table for the matching rates shows that the theoretical matching rates calculated by the algorithm are quite close to the simulated (actual) matching rates, already for moderate values of $\lambda$. The results for the waiting times and idle times confirm our intuition that they should converge to the targeted quality of service requirements: for large values of $\lambda$, the probability mass of the waiting times in ED concentrates near $W$ and the probability mass of the idle times in QD concentrates near $T$. In the QED regime, waiting times, idle times and abandonment rates are small.

\subsection{ Example 2 -- \boldmath $5\times5$ Hamiltonian  Graph with Differentiated Service}
\label{sec.example2}

This example illustrates differentiated service, with the customer types divided into three classes: high, standard and low priority customers.

\begin{center}
\begin{tabular}{|p{0.92\textwidth}|}
\hline
\hfill\textbf{Example 2 -- System and Data}\hfill\mbox{}\\
\hline
There are 5 types of customers and 5 types of servers. The total arrival rate is parameterized by $\lambda$.  The graph and the values of $\alpha_{c_i}\lambda$ are described in the following figure:

   \begin{center}
    \includegraphics[scale=0.30]{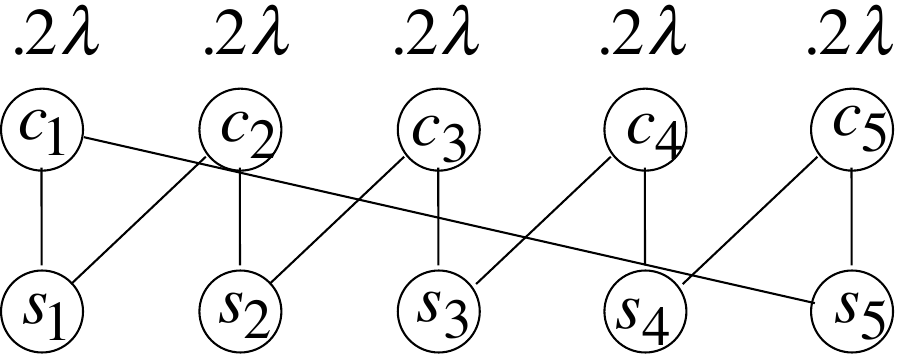}
    \end{center}

The fractions $\alpha_{c_i}$ are all equal, i.e., $\alpha_{c_i} = \frac15$.
The patience times are all exponentially distributed with mean 10.
The service times are all uniformly distributed, with parameters as given in the table below.

\begin{center}
\begin{tabular}{|c|ccccc|}
\hline
\multicolumn{6}{|c|}{Service time distributions}\\
\hline
$G_{c_i,s_j}$    & $c_1$ & $c_2$ & $c_3$ & $c_4$ & $c_5$ \\
\hline
$s_1$ & U(2, 6) & U(2, 4) & & & \\
$s_2$ & & U(1, 3) & U(4, 7) & &  \\
$s_3$ & & & U(3, 6) & U(2, 6) & \\
$s_4$ & & & & U(1, 5) & U(6, 11) \\
$s_5$ & U(3, 7) & & &  & U(4, 9) \\
\hline
\end{tabular}
\quad
\begin{tabular}{|c|ccccc|}
\hline
\multicolumn{6}{|c|}{Mean service times}\\
\hline
$G_{c_i,s_j}$    & $c_1$ & $c_2$ & $c_3$ & $c_4$ & $c_5$ \\
\hline
$s_1$ & 8 & 3 & & & \\
$s_2$ & & 2 & 5.5 & &  \\
$s_3$ & & & 4.5 & 4 & \\
$s_4$ & & & & 3 & 8.5 \\
$s_5$ & 5 & & &  & 6.5 \\
\hline
\end{tabular}
\end{center}

Only the \emph{mean} service times are used by the design algorithms. The full distributions are used in the simulations.\\
\hline
\end{tabular}
\end{center}

\begin{figure}[htb]
   \begin{center}
    \includegraphics[scale=0.30]{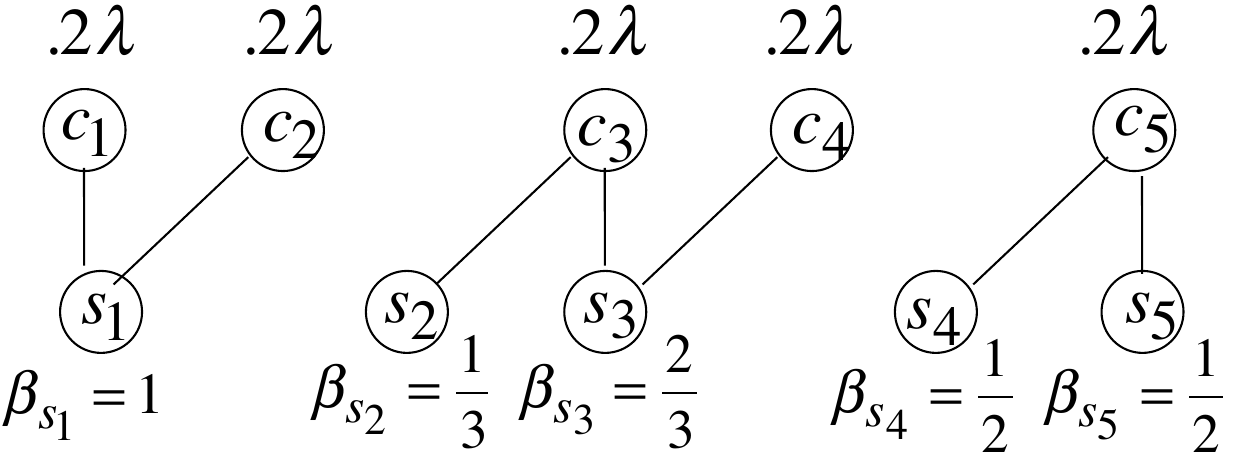}
    \caption{Decomposed system of the system of Example 2.}\label{fig:dec}
   \end{center}
\end{figure}

We consider  the following  decomposition of the system of Example 2:
\[
\C^{(1)} = \{c_1,c_2\},\;\; \S^{(1)} =\{s_1\},\quad
\C^{(2)} = \{c_3,c_4\},\;\; \S^{(2)} =\{s_2,s_3\},\quad
\C^{(3)} = \{c_5\},\;\; \S^{(3)} =\{s_4,s_5\}.
\]
The decomposed system is described in Figure \ref{fig:dec}.
Note that this decomposition results from eliminating three links in the compatibity graph: the link from $s_2$ to $c_2$, $s_4$ to $c_4$, and from $s_5$ to $c_1$.
We then have that:\\

\begin{compactitem}
\item[-]
$\C^{(1)},\S^{(1)}$ in isolation is a ``V'' system, with arrival rate $0.4 \lambda$, adjusted fractions $\alpha_{c_1}=\alpha_{c_2} = \frac12$ and $\beta_{s_1} = 1$.
\item[-]
$\C^{(2)},\S^{(2)}$ in isolation is an ``N'' system, with arrival rates $0.4\lambda$, adjusted fractions $\alpha_{c_3}=\alpha_{c_4} = \frac12$ and we take $\beta_{s_2} = \frac13, \beta_{s_3}=\frac23$.
\item[-]
$\C^{(3)},\S^{(3)}$ in isolation is a ``$\Lambda$'' system, with arrival rate  $\lambda=0.2 \lambda$, $\alpha_{c_5} = 1$
and we take $\beta_{s_4} = \beta_{s_5}=\frac12$.
\end{compactitem}

We note that the subnetworks have simple compatibility structure so that matching rates are obtained immediately as:
\[
r_{c_1,s_1}^{(1)}=r_{c_2,s_1}^{(1)}=\frac{1}{2},\quad
r_{c_3,s_2}^{(2)}=r_{c_3,s_3}^{(2)}=r_{c_4,s_3}^{(2)}=\frac{1}{3},\quad
r_{c_5,s_4}^{(3)}=r_{c_5,s_5}^{(3)}=\frac{1}{2}.
\]

We make three designs for this network, in which customers $c_1,c_2$ have high priority, $c_3,c_4$ have standard priority, and $c_5$ have low priority.   The first design is for a system in ED regime, the second is for a mixed design with the top priority sub-system in QD, the middle priority sub-system in QED and low priority sub-system in ED, and the third design is for the three systems in QD regime. It is readily verified, by checking Conditions (\ref{eqn.crp}), that in each design complete resource pooling holds for each subsystem.
Table \ref{tbl:wf2} shows the calculated workforce required for each type of server for the three designs, as a function of $\lambda$.

\begin{table}
\begin{center}
\begin{tabular}{|c|ccccc|ccccc|ccccc|}
\hline
 & \multicolumn{15}{|c|}{Required workforce}\\
 \hline
& \multicolumn{5}{|c|}{ED regime} & \multicolumn{5}{|c|}{Mixed regime} & \multicolumn{5}{|c|}{QD regime} \\
&   $W=1$ & \multicolumn{2}{c}{$W=2$} & \multicolumn{2}{c|}{$W=3$}
&   $T=1$ & \multicolumn{2}{c}{QED} & \multicolumn{2}{c|}{$W=1$}
&   $T=2$ & \multicolumn{2}{c}{$T=1$} & \multicolumn{2}{c|}{$T=0.5$}  \\
\hline
$\lambda$
  & $n_{s_1}$ & $n_{s_2}$ & $n_{s_3}$ & $n_{s_4}$ & $n_{s_5}$
  & $n_{s_1}$ & $n_{s_2}$ & $n_{s_3}$ & $n_{s_4}$ & $n_{s_5}$
  & $n_{s_1}$ & $n_{s_2}$ & $n_{s_3}$ & $n_{s_4}$ & $n_{s_5}$
  \\
\hline
20 &  25 & 12 & 18 & 13 & 10 & 36 & 15 & 22 & 15 & 12 & 44 & 17 & 27 & 18 & 14 \\
40 &  51 &  24 & 36 & 25 & 19 & 72 & 29 & 44 & 31 & 24 & 88 & 35 & 55 & 36 & 28 \\
60 &  76 & 36 &   54 & 38 & 29 & 108 & 44 & 66 & 46 & 35 & 132 & 52 & 82 & 54 & 42 \\
100 & 127 & 60 & 90 &  63 & 48 & 180 & 73 & 110 & 77 & 59 & 220 & 87 & 137 & 90 & 70 \\
200 & 253 & 120 & 180 & 126 & 96 & 360 & 147 & 220 & 154 & 118 & 440 & 173 & 273 & 180 & 140 \\
\hline
\end{tabular}
\caption{Calculated required workforce for Example 2}\label{tbl:wf2}
\end{center}
\end{table}

The simulation results for Example 2 are listed in the tables below. The results illustrate that the system with differentiated service under FCFS-ALIS performs in accordance with the design.

\begin{center}
\small
\setlength{\tabcolsep}{1mm}
\begin{tabular}{|c|ccccc|ccccc|ccccc|}
\hline
 & \multicolumn{15}{|c|}{Matching rates}\\
\hline
& \multicolumn{5}{|c|}{ED regime} & \multicolumn{5}{|c|}{Mixed regime} & \multicolumn{5}{|c|}{QD regime} \\
\hline
 & \multicolumn{15}{|c|}{Theoretical}\\
\hline
$r_{c_i,s_j}$    & $c_1$ & $c_2$ & $c_3$ & $c_4$ & $c_5$ & $c_1$ & $c_2$ & $c_3$ & $c_4$ & $c_5$& $c_1$ & $c_2$ & $c_3$ & $c_4$ & $c_5$\\
\hline
$s_1$ &   0.216 & 0.216 &       &       &       &  0.204 & 0.204 &       &       &       &  0.200 & 0.200 &       &       &       \\
$s_2$ &         &       & 0.130 &       &       &        &       & 0.136 &       &       &        &       & 0.133 &       &       \\
$s_3$ &         &       & 0.065 & 0.195 &       &        &       & 0.068 & 0.204 &       &        &       & 0.067 & 0.200 &       \\
$s_4$ &         &       &       &       & 0.088 &        &       &       &       & 0.092 &        &       &       &       & 0.100 \\
$s_5$ &         &       &       &       & 0.088 &        &       &       &       & 0.092 &        &       &       &       & 0.100 \\
\hline
 & \multicolumn{15}{|c|}{$\lambda=20$}\\
\hline
$r_{c_i,s_j}$    & $c_1$ & $c_2$ & $c_3$ & $c_4$ & $c_5$ & $c_1$ & $c_2$ & $c_3$ & $c_4$ & $c_5$& $c_1$ & $c_2$ & $c_3$ & $c_4$ & $c_5$\\
\hline
$s_1$ &   0.215 & 0.196 &       &       &       &  0.207 & 0.194 &       &       &       &  0.200 & 0.186 &       &       &       \\
$s_2$ &         & 0.021 & 0.121 &       &       &        & 0.014 & 0.126 &       &       &        & 0.015 & 0.118 &       &       \\
$s_3$ &         &       & 0.079 & 0.176 &       &        &       & 0.077 & 0.186 &       &        &       & 0.083 & 0.181 &       \\
$s_4$ &         &       &       & 0.017 & 0.084 &        &       &       & 0.014 & 0.086 &        &       &       & 0.020 & 0.094 \\
$s_5$ &   0.002 &       &       &       & 0.089 &  0.001 &       &       &       & 0.095 &  0.002 &       &       &       & 0.100 \\
\hline
 & \multicolumn{15}{|c|}{$\lambda=200$}\\
\hline
$r_{c_i,s_j}$    & $c_1$ & $c_2$ & $c_3$ & $c_4$ & $c_5$ & $c_1$ & $c_2$ & $c_3$ & $c_4$ & $c_5$& $c_1$ & $c_2$ & $c_3$ & $c_4$ & $c_5$\\
\hline
$s_1$ &   0.216 & 0.216 &       &       &       & 0.205 & 0.205 &       &       &       &  0.200 & 0.200 &       &       &       \\
$s_2$ &         &       & 0.130 &       &       &       &       & 0.135 &       &       &        &       & 0.132 &       &       \\
$s_3$ &         &       & 0.066 & 0.194 &       &       &       & 0.068 & 0.201 &       &        &       & 0.069 & 0.197 &       \\
$s_4$ &         &       &       & 0.001 & 0.088 &       &       &       & 0.001 & 0.093 &        &       &       & 0.003 & 0.099 \\
$s_5$ &         &       &       &       & 0.088 &       &       &       &       & 0.093 &        &       &       &       & 0.101 \\
\hline
\end{tabular}
\end{center}

\begin{center}
\begin{tabular}{|c|c|c|c|}
\hline
 & \multicolumn{3}{|c|}{Customer waiting times and server idle times}\\
\hline
ED Regime & Mixed Regime & Mixed Regime & QD Regime\\
Waiting times & Waiting times & Idle times & Idle times\\
\hline
&&&\\[-2mm]
\includegraphics[width=0.24\textwidth]{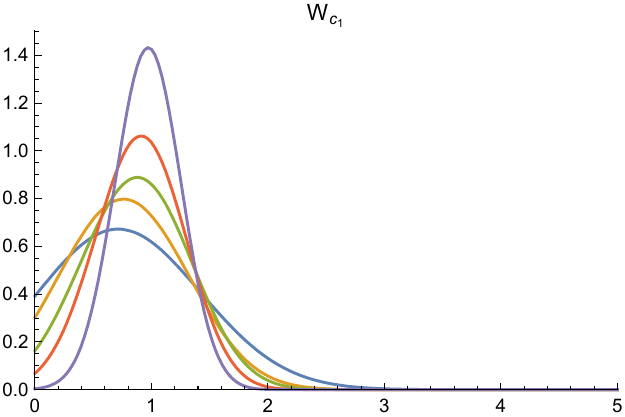} &                                                         &\includegraphics[width=0.24\textwidth]{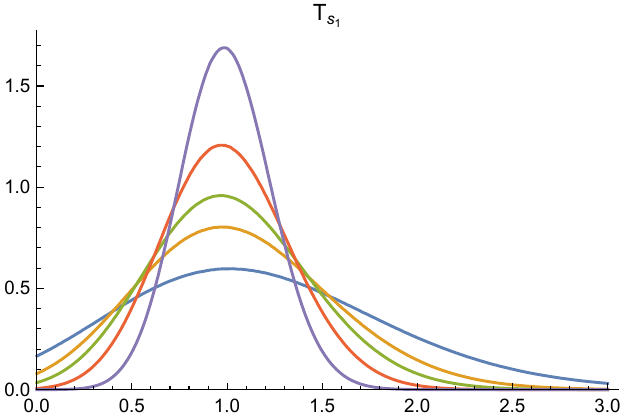} &\includegraphics[width=0.24\textwidth]{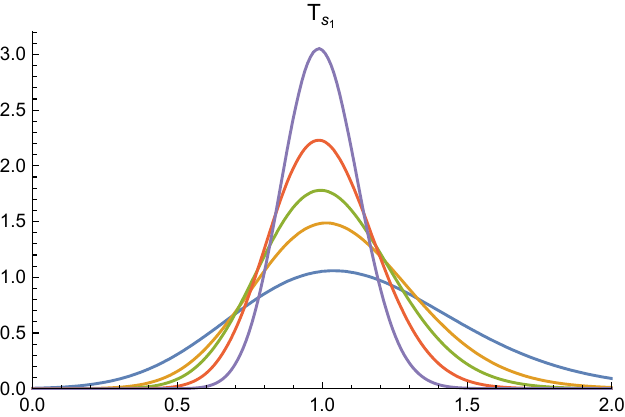} \\
\includegraphics[width=0.24\textwidth]{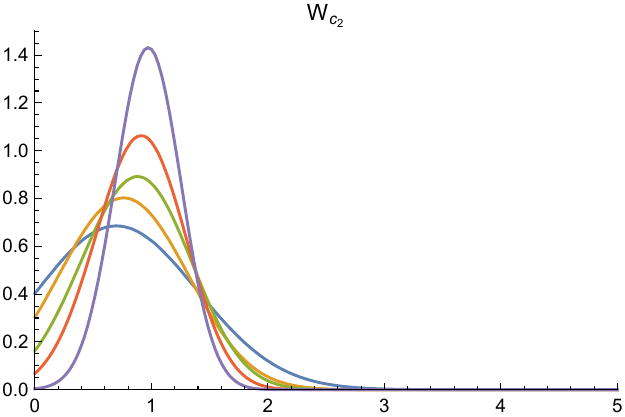} &                                                         &\includegraphics[width=0.24\textwidth]{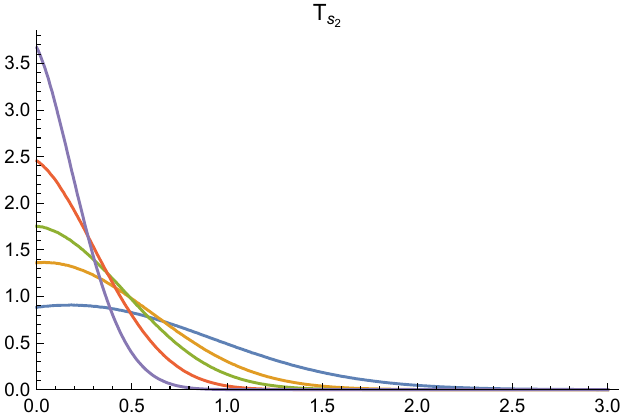} &\includegraphics[width=0.24\textwidth]{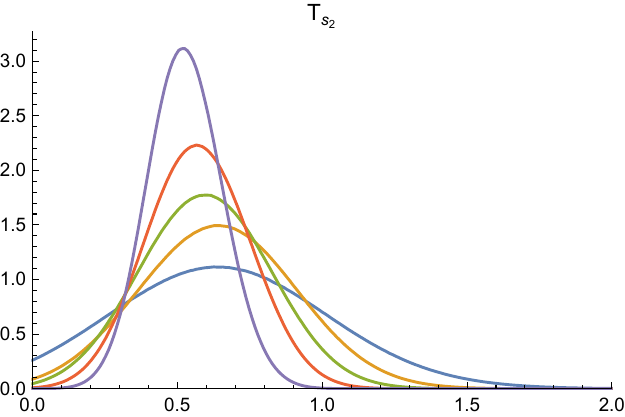}   \\
\includegraphics[width=0.24\textwidth]{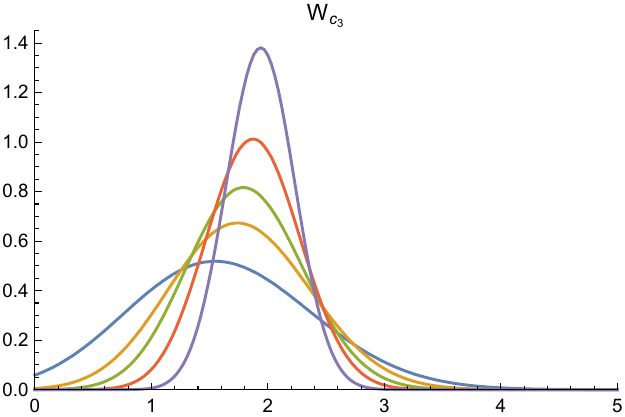} & \includegraphics[width=0.24\textwidth]{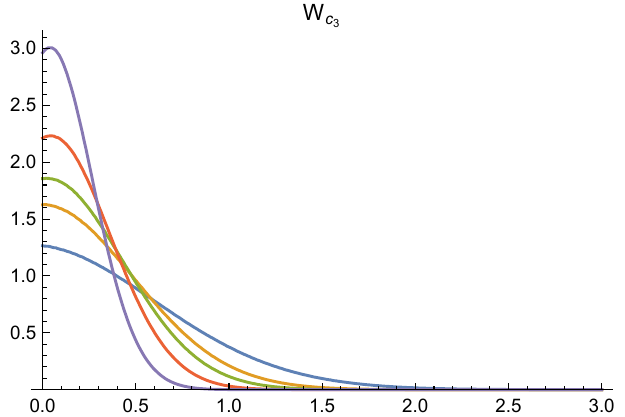}  &\includegraphics[width=0.24\textwidth]{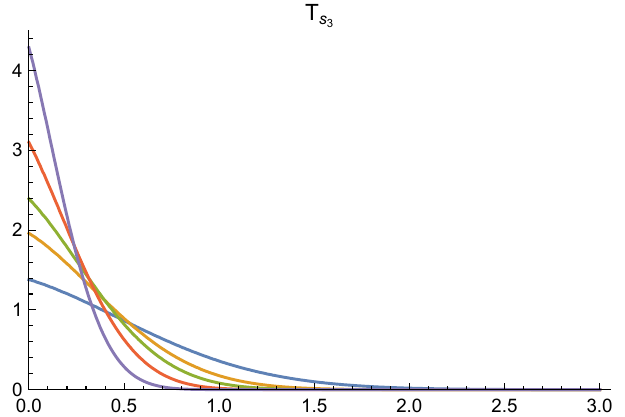} &\includegraphics[width=0.24\textwidth]{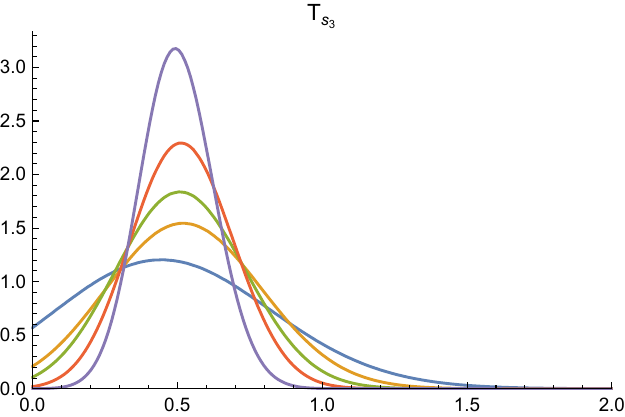}     \\
\includegraphics[width=0.24\textwidth]{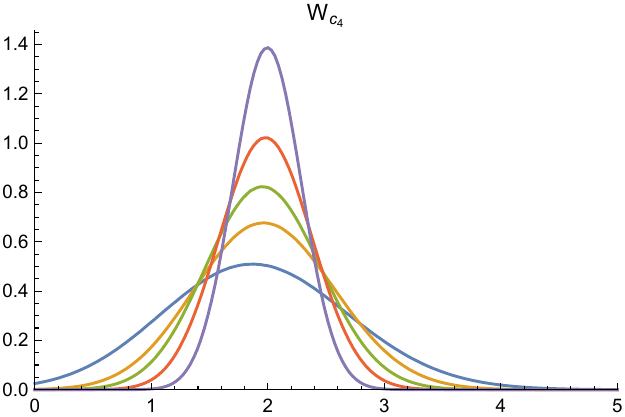} & \includegraphics[width=0.24\textwidth]{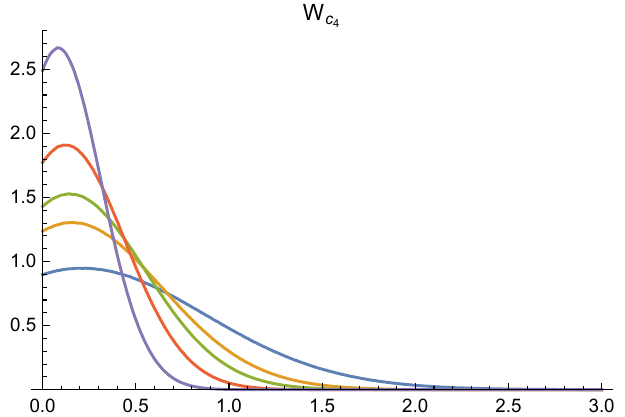}  &                                                    &\includegraphics[width=0.24\textwidth]{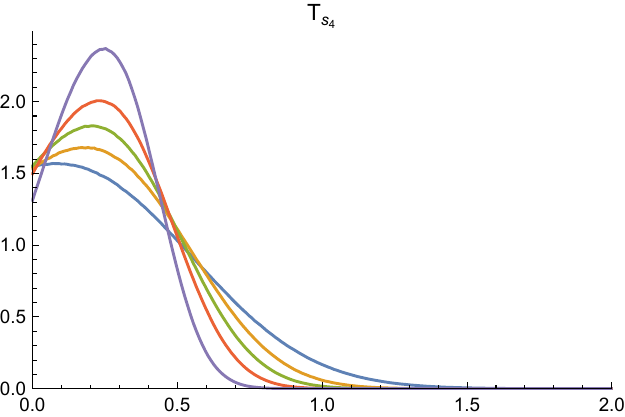}       \\
\includegraphics[width=0.24\textwidth]{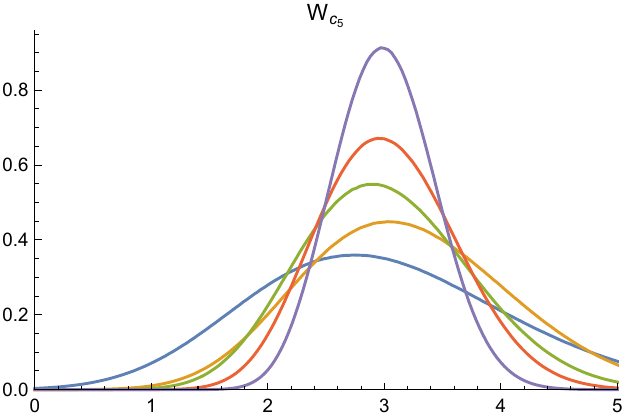} & \includegraphics[width=0.24\textwidth]{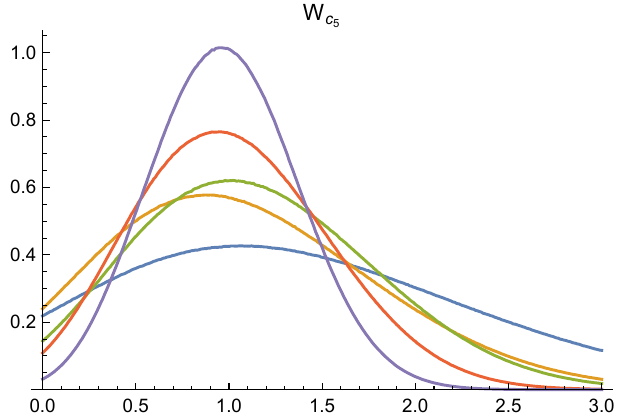}  &                                                    &\includegraphics[width=0.24\textwidth]{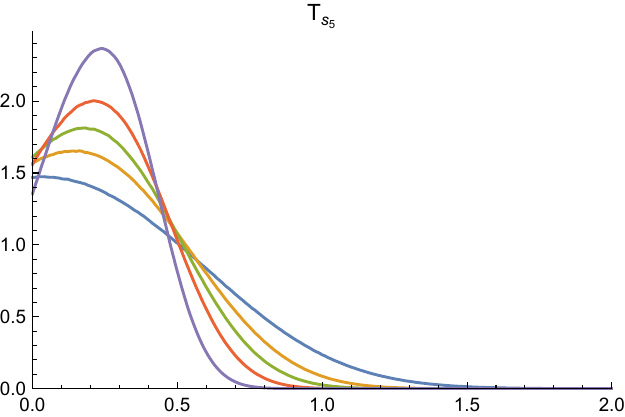}         \\
\hline
\multicolumn{4}{|c|}{\includegraphics[width=0.6\textwidth]{legend2.pdf}}\\
\hline
\end{tabular}
\end{center}
%
%
%


\begin{table}[h]
\begin{center}
\begin{tabular}{|l|c|c|c|c|c|c|}
\hline
 & \multicolumn{6}{|c|}{Fraction of no wait and of no idling}\\
\hline
    & \multicolumn{2}{|c|}{ED regime}&    \multicolumn{2}{|c|}{Mixed regime}& \multicolumn{2}{|c|}{QD regime} \\
\hline
$\lambda$    & No waiting & No idling & No waiting &  No idling  & No waiting   & No idling \\
\hline
20   &  0.056 & 0.934  &  0.578 & 0.398 &     0.869 & 0.124 \\
60   &  0.010 & 0.988  &  0.570 & 0.412 &     0.935 & 0.063 \\
200  &  0.000 & 1.000  &  0.556 & 0.430 &     0.974 & 0.025 \\
\hline
\end{tabular}
\end{center}
\caption{Simulation results for Example 2}\label{tabl:res2}
\end{table}

\subsection{Example 3 -- \boldmath$6\times6$ Symmetric Degree 3 Graph with Pooled Service}
\label{subsect:example3}

We now consider a more complex graph to examine the validity of the matching rates conjecture.

\label{sec.example3}
\begin{center}
\begin{tabular}{|p{0.92\textwidth}|}
\hline
\hfill\textbf{Example 3 -- System and Data}\hfill\mbox{}\\
\hline
There are 6 types of customers and 6 types of servers. The total arrival rate is parameterized by $\lambda$,  The graph and the values of $\lambda_{c_i}/\lambda$ are described in the following figure:

\parbox{0.45\textwidth}{
   \begin{center}
    \includegraphics[scale=0.30]{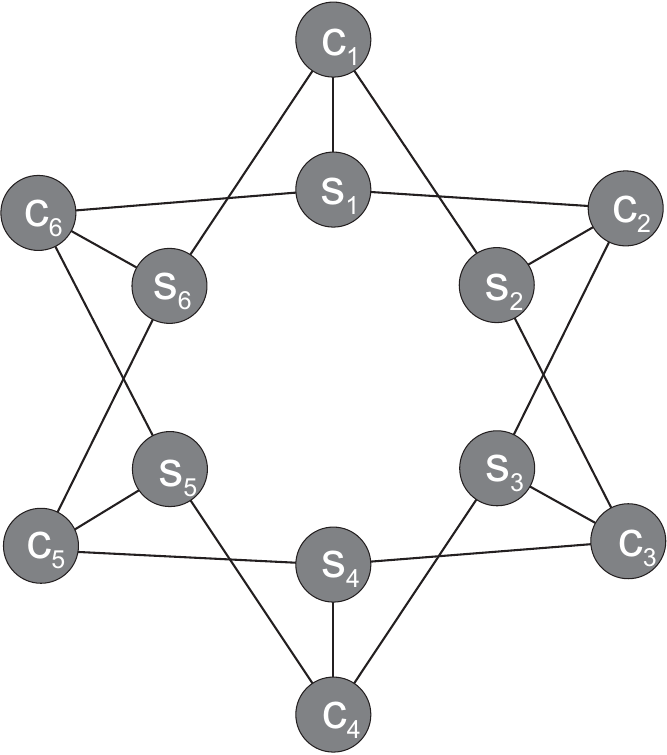}
    \end{center}
}
\parbox{0.45\textwidth}{
$$
\begin{array}{|l|l|}
\hline
\multicolumn{2}{|c|}{\lambda_{c_i}/\lambda}\\
\hline
c_1 & 1/9 \\
c_2 & 2/9 \\
c_3 & 1/9 \\
c_4 & 2/9 \\
c_5 & 1/9 \\
c_6 & 2/9 \\
\hline
\end{array}
$$
}

The patience times are all exponentially distributed with mean 10.
The service times distributions are given in the table below.

\hspace{1.5cm}
\parbox{0.45\textwidth}{
\begin{tabular}{|c|cccccc|}
\hline
\multicolumn{7}{|c|}{Service time distributions}\\
\hline
    & $c_1$ & $c_2$ & $c_3$  & $c_4$ & $c_5$ & $c_6$\\
\hline
$s_1$ & U & E &   &   &   & P \\
$s_2$ & P & U & E &   &   &   \\ $s_3$ &   & P & U & E &   &   \\ $s_4$ &   &   & P & U & E &   \\ $s_5$ &   &   &   & P & U & E \\ $s_6$ & E &   &   &   & P & U \\ \hline
\end{tabular}
}
\parbox{0.35\textwidth}{
With:\\
$E\sim \textrm{Exp}(1/4)$, mean 4\\
$P\sim \textrm{Pareto}(3, 3)$, mean 4.5\\
$U\sim \textrm{U}(1, 3)$, mean 2 \\
}
Only the \emph{mean} service times are used by the design algorithms. The full distributions are used in the simulations.\\
\hline
\end{tabular}
\end{center}

In the designs for Example 3 we take as service fractions:  $\beta_{s_j}=1/6, \quad j=1,\ldots,6$. It then follows, by checking Conditions (\ref{eqn.crp}), that in each regime (ED, QED and QD) service is pooled. We take $W=1$ in the ED regime and $T=0.5$ in the QD regime.

\begin{table}[htb]
\begin{center}
\setlength{\tabcolsep}{1mm}
\begin{tabular}{|c|cccccc|cccccc|cccccc|}
\hline
 & \multicolumn{18}{|c|}{Required workforce}\\
 \hline
& \multicolumn{6}{|c|}{ED regime} & \multicolumn{6}{|c|}{QED regime} & \multicolumn{6}{|c|}{QD regime} \\
\hline
$\lambda$
  & $n_{s_1}$ & $n_{s_2}$ & $n_{s_3}$ & $n_{s_4}$ & $n_{s_5}$& $n_{s_6}$
  & $n_{s_1}$ & $n_{s_2}$ & $n_{s_3}$ & $n_{s_4}$ & $n_{s_5}$& $n_{s_6}$
  & $n_{s_1}$ & $n_{s_2}$ & $n_{s_3}$ & $n_{s_4}$ & $n_{s_5}$& $n_{s_6}$
  \\
\hline
20  &  12 & 9 & 12 & 9 & 12 & 9 &         13 & 10 & 13 & 10 & 13 & 10 &         15 & 12 & 15 & 12 & 15 & 12 \\
40  &  23 & 19 & 23 & 19 & 23 & 19 &      26 & 21 & 26 & 21 & 26 & 21 &         29 & 24 & 29 & 24 & 29 & 24 \\
60  &  35 & 28 & 35 & 28 & 35 & 28 &      39 & 31 & 39 & 31 & 39 & 31 &         44 & 36 & 44 & 36 & 44 & 36 \\
100 &  58 & 47 & 58 & 47 & 58 & 47 &      64 & 52 & 64 & 52 & 64 & 52 &         73 & 60 & 73 & 60 & 73 & 60 \\
200 &  117 & 94 & 117 & 94 & 117 & 94 &   129 & 104 & 129 & 104 & 129 & 104 &   146 & 121 & 146 & 121 & 146 & 121 \\
\hline
\end{tabular}
\caption{Calculated workforce for Example 3}\label{tbl:wf3}
\end{center}
\end{table}

The theoretical matching rates are given in the table below. Note that these matching rates are the same for all three regimes (ED, QED and QD). This is caused by the fact that all customers have the same patience distribution and the same target waiting times, resulting in equal $\alpha_{c_i}$ for all three regimes.
The simulated matching rates are also practically identical for all three regimes. In the tables we depict the averages over the three regimes, but the actual differences between the three simulated values and their averaged values are less than 0.001.

\begin{table}{htb}
\begin{center}
\begin{tabular}{|c|cccccc|}
\hline
 \multicolumn{7}{|c|}{Theoretical matching rates (ED, QED, QD)} \\
\hline
$r_{c_i,s_j}$    & $c_1$ & $c_2$ & $c_3$ & $c_4$ & $c_5$  & $c_6$\\
\hline
$s_1$ &  0.028 & 0.069 &       &       &       & 0.069  \\
$s_2$ &  0.041 & 0.084 & 0.041 &       &       &        \\
$s_3$ &        & 0.069 & 0.028 & 0.069 &       &        \\
$s_4$ &        &       & 0.041 & 0.084 & 0.041 &        \\
$s_5$ &        &       &       & 0.069 & 0.028 & 0.069  \\
$s_6$ &  0.041 &       &       &       & 0.041 & 0.084  \\
\hline
 \multicolumn{7}{|c|}{Simulated matching rates $\lambda=20$ (Average ED, QED, QD)} \\
\hline
$r_{c_i,s_j}$    & $c_1$ & $c_2$ & $c_3$ & $c_4$ & $c_5$  & $c_6$\\
\hline
$s_1$ &   0.030 & 0.071 &       &       &       & 0.070 \\
$s_2$ &   0.041 & 0.080 & 0.041 &       &       &       \\
$s_3$ &         & 0.070 & 0.030 & 0.071 &       &       \\
$s_4$ &         &       & 0.041 & 0.080 & 0.041 &       \\
$s_5$ &         &       &       & 0.070 & 0.030 & 0.071 \\
$s_6$ &   0.041 &       &       &       & 0.041 & 0.080 \\
\hline
 \multicolumn{7}{|c|}{Simulated matching rates $\lambda=200$ (Average ED, QED, QD)} \\
\hline
$r_{c_i,s_j}$    & $c_1$ & $c_2$ & $c_3$ & $c_4$ & $c_5$  & $c_6$\\
\hline
$s_1$ &    0.028 & 0.069 &       &       &       & 0.069 \\
$s_2$ &    0.041 & 0.084 & 0.041 &       &       &       \\
$s_3$ &          & 0.069 & 0.028 & 0.069 &       &       \\
$s_4$ &          &       & 0.041 & 0.084 & 0.041 &       \\
$s_5$ &          &       &       & 0.069 & 0.028 & 0.069 \\
$s_6$ &    0.041 &       &       &       & 0.041 & 0.084 \\
\hline
\end{tabular}
\end{center}

%
%
%

\begin{center}
\begin{tabular}{|l|c|c|c|c|c|c|}
\hline
 & \multicolumn{6}{|c|}{Fraction of no wait and of no idling} \\
\hline
    & \multicolumn{2}{|c|}{ED regime}&    \multicolumn{2}{|c|}{QED regime}& \multicolumn{2}{|c|}{QD regime} \\
\hline
$\lambda$    & No waiting & No idling & No waiting &  No idling  & No waiting   & No idling \\
\hline
20   &     0.044 & 0.950 &       0.278 & 0.709 &  0.862 & 0.135 \\
40   &     0.014 & 0.985 &       0.398 & 0.594 &  0.929 & 0.071 \\
60   &     0.003 & 0.997 &       0.351 & 0.642 &  0.974 & 0.026 \\
100  &     0.000 & 1.000 &       0.314 & 0.681 &  0.994 & 0.006 \\
200  &     0.000 & 1.000 &       0.352 & 0.644 &  1.000 & 0.000 \\
\hline
\end{tabular}
\end{center}
\caption{Simulation results for Example 3}\label{tbl:res3}
\end{table}

Finally, we depict the abandonment rates in the Table \ref{tbl:ab3}. Due to the symmetry in the system, we only report abandonment rates for customer types $c_1$ and $c_2$, because the results for customer types 3 and 5 are identical to those of customer type 1, and those of customer types 4 and 6 are identical to those of type 2.

\begin{table}[htb]
\begin{center}
\begin{tabular}{|l|c|c|c|c|c|c|}
\hline
 & \multicolumn{6}{|c|}{Simulated abandonment rates} \\
\hline
    & \multicolumn{2}{|c|}{ED regime}&    \multicolumn{2}{|c|}{QED regime}& \multicolumn{2}{|c|}{QD regime} \\
\hline
$\lambda$ & $c_1$ & $c_2$ & $c_1$ & $c_2$ & $c_1$ & $c_2$ \\
\hline
20 &  0.104 & 0.106 &  0.041 & 0.043 &     0.003 & 0.003 \\
200 & 0.094 & 0.095 &  0.010 & 0.011 &     0.000 & 0.000 \\
\hline
\end{tabular}
\end{center}
\caption{Simulated abandonment rates for Example 3}\label{tbl:ab3}
\end{table}

The expected abandonment rates are 0.095 for every customer type (due to symmetry) in the ED regime and zero in the QED and QD regimes.  The pattern here is as expected: for $\lambda=20$ the simulated abandonment rates are quite close to the desired values, in particular in the ED and QD regimes. For $\lambda=200$ we are also very close to achieving our target abandonment rates in the QED regime. In this example we have taken all patience times to be exponentially distributed with mean 10, which is quite large compared to the target waiting times (no waiting in the QD and QED regimes, and $W=1$ in the ED regime). In the next sub-section we conduct a more in-depth study of the impact of the patience time distributions on the abandonment rates.

\subsection{Impact of the patience-time distribution.}
\label{subsect:differentpatiencedists}

So far we have not studied the impact of the \emph{distribution} of the customer patience time. In order to gain more insight in this topic, we take the setting of Example 3 because of its symmetry, which makes it suitable for measuring the impact of the patience-time distribution on the various performance measures. Customers of types 1 and 2 have exponentially distributed patience times; types 3 and 4 have a uniform distribution; types 5 and 6 have a Pareto distribution. The means of these patience times are 2.5 for customer types 1, 3, 5 and 3.0 for types 2, 4 and 6.

\begin{table}[ht]
\begin{center}
\begin{tabular}{|c|c|c|}
\hline
\multicolumn{2}{|c|}{Patience time distributions} & Expected abandonment rates ED  \\
\hline
& $F_{c_i}$ & $F_{c_i}(W)$   \\
\hline
$c_1$ & Exp(0.4)     & 0.330\\
$c_2$ & Exp(0.333) & 0.283\\
$c_3$ & U(0.5, 4.5) & 0.125\\
$c_4$ & U(1.0, 5.0) & 0.000\\
$c_5$ & Pareto(0.8333, 1.5)& 0.239 \\
$c_6$ & Pareto(1, 1.5) & 0.000\\
\hline
\end{tabular}
\end{center}
\caption{Patience time distributions and expected abandonment rates in the ED regime (with $W=1$) for Example 3}
\label{tbl:differentpatiencedistsinput}
\end{table}

Table \ref{tbl:differentpatiencedistsinput} gives a detailed overview of the distributions and the parameters. We also show the expected fractions of abandonments in the ED regime (with $W=1$). Recall that these expected abandonment rates are zero in the QED regime. For this reason, the required numbers of servers in the QED regime do not depend on the patience time distribution, implying they are identical to the numbers listed in Table \ref{tbl:wf3}. The staffing levels in the ED regime are slightly different due to the different distributions. For $\lambda=20$ these differences are still very small (difference of at most one server), but for $\lambda=200$ the number of servers of each type are typically smaller than in the previous sub-section:
\[
n_{s_1}=115, n_{s_2}=91, n_{s_3}=110, n_{s_4}=85, n_{s_5}=113, n_{s_6}=84.
\]
This can be explained by the fact that the abandonment rates for most customers types are much higher in this example than in the previous example. We have deliberately chosen higher values to emphasize the impact of the patience-time distributions.

We ran the same number of simulations (1,000 $\times$ 1,000,000 matches) as before. The simulation results (which we have omitted for reasons of compactness) clearly indicate that the different patience time distributions do \emph{not} have an impact on the simulated matching rates nor or the mean waiting times, which are still extremely close to the theoretical values -- even for small $\lambda$. However, Table \ref{tbl:abandonmentratesdifferentpatiencedist}, with simulated fractions of abandonments, shows some interesting patterns.
First, we observe that the simulated arrival rates clearly converge towards their target values. Interestingly, this convergence is slower in those cases where we wish to achieve zero abandonments in the ED regime. The ED regime is not primarily designed to favor customers and a target of 0\% abandonments is difficult to achieve in this regime. The second interesting observation, which is in contrast to the previous remark, is that the uniform and Pareto distribution make it relatively easy to achieve 0\% abandonments in the QED regime. The reason for this phenomenon is that the support of both distributions has a  positive offset. As $\lambda$ increases, the tail of the waiting time distribution becomes so light, that no customer experiences a waiting time that is longer than the smallest possible patience. All in all we can conclude that:
\begin{itemize}
\item Customer waiting times and server idle times are (nearly) insensitive to the patience time distributions;
\item Abandonment rates are close to their target values irrespective of the patience time distributions, in particular for large $n$;
\item Abandonment rates of 0\% are difficult to achieve in the ED regime, but easy to achieve in the QED regime if the support of the patience time distribution has a positive offset. For completeness we stress that the QD regime always results in zero abandonments, even for small $\lambda$.
\end{itemize}

\begin{table}[ht]
\begin{tabular}{|c|cccccc|cccccc|}
\hline
 & \multicolumn{12}{|c|}{Simulated abandonment rates} \\
\hline
    & \multicolumn{6}{|c|}{ED regime}&    \multicolumn{6}{|c|}{QED regime} \\
\hline
$\lambda$  & $c_1$ & $c_2$ & $c_3$  & $c_4$ & $c_5$ & $c_6$ & $c_1$ & $c_2$ & $c_3$  & $c_4$ & $c_5$ & $c_6$ \\
\hline
 20  &   0.280 & 0.246 & 0.112 & 0.046 & 0.238 & 0.133 & 0.114 & 0.101 & 0.023 & 0.004 & 0.031 & 0.013 \\
 40  &   0.291 & 0.253 & 0.104 & 0.025 & 0.198 & 0.088 & 0.065 & 0.058 & 0.005 & 0.000 & 0.004 & 0.001 \\
 60  &   0.309 & 0.268 & 0.113 & 0.022 & 0.216 & 0.085 & 0.064 & 0.055 & 0.003 & 0.000 & 0.001 & 0.000 \\
 100  &  0.317 & 0.274 & 0.117 & 0.016 & 0.217 & 0.069 & 0.057 & 0.049 & 0.001 & 0.000 & 0.000 & 0.000 \\
 200  &  0.325 & 0.280 & 0.122 & 0.011 & 0.226 & 0.054 & 0.036 & 0.031 & 0.000 & 0.000 & 0.000 & 0.000 \\
\hline
\end{tabular}
\caption{Abandonment rates for the model with different patience distributions in Section~\ref{subsect:differentpatiencedists}.}
\label{tbl:abandonmentratesdifferentpatiencedist}
\end{table}

\subsection{Example calculations for given $n_{s_j}/n$}
\label{sec.minimization}

As discussed before, there is an alternative way to determine the staffing levels, which does not require an explicit choice of the $\beta_{s_j}$'s. In some practical cases it may be more natural to specify the desired fractions of the total number of servers of each type instead, denoted by
\[\underline{\theta}_s=(\theta_{s_1},\theta_{s_2},\dots,\theta_{s_I})=(n_{s_1}, \dots, n_{s_I})/n.\]
In this section we consider staffing decisions made on the basis of specification of the $\theta_{s_j}$'s. We revisit Example~1 from Section \ref{sec.example1} and Example~3 from Section \ref{sec.example3}.  We note that in Example~2 from Section \ref{sec.example2}, the one to one relation between $\beta_{s_j}$ and $\theta_{s_j}$ is immediate and does not require numerical minimization.

\paragraph{Example 1 revisited.} In this example we consider the same system with three customer types and three server types. In Section \ref{sec.example1} we have chosen settings that ensure complete resource pooling, namely $(\beta_{s_1}, \beta_{s_2}, \beta_{s_3})=(0.3, 0.3, 0.4)$. Instead of specifying $(\beta_{s_1}, \beta_{s_2}, \beta_{s_3})$ we now specify desired fractions of server types $(\theta_{s_1}, \theta_{s_2}, \theta_{s_3})$, distinguishing between three cases:
\begin{itemize}
\item $\theta_{s_1}=\theta_{s_2}=\theta_{s_3}=1/3$, i.e. all fractions are equal;
\item $\theta_{s}=(1/6,1/3,1/2)$, i.e. half of the servers should be of type 3, one third of type 2, and the rest of type 1;
\item $\theta_{s}=(1/2,1/3,1/6)$, i.e. half of the servers should be of type 1, one third of type 2, and the rest of type 3.
\end{itemize}
In order to find a vector $\beta_{s}$ that results in the desired $\theta_s$, we numerically minimize the function $\Delta$ as defined in Equation \eqref{eqn:Delta}. When minimizing $\Delta$, we impose additional restrictions \eqref{eqn.crp} to ensure complete resource pooling. For this simple network, these restrictions can be written out as follows:
\begin{eqnarray*}
\mbox{QD, QED:} && 0 < \beta_{s_1} < 0.7, \qquad  0.5 < \beta_{s_1}+\beta_{s_2} < 1, \qquad 0 < \beta_{s_2} < 0.8; \\
\mbox{ED:} &&      0 < \beta_{s_1} < 0.720, \quad 0.513 < \beta_{s_1}+\beta_{s_2} < 1, \quad 0 < \beta_{s_2} < 0.794.
\end{eqnarray*}
The results of the numerical routine, for all three regimes, can be found in Table \ref{tbl:example1thetas}. The first column gives the vector $\theta_s$ of desired fractions of servers, the second column indicates whether this desired fraction could actually be attained. The third column gives the vector $\beta_s$ that minimizes Equation \eqref{eqn:Delta}, and the last column displays the required staffing levels of each server type when $\lambda=100$.

\begin{table}[htb]
\begin{tabular}{|c|c|c|c|}
\hline
Desired $\theta_{s}$ & Realized & $\beta_s$ & Staffing levels for $\lambda=100$ \\
\hline
\multicolumn{4}{|c|}{\bfseries ED Regime}\\
\hline
(1/3, 1/3, 1/3) & Yes & (0.213803, 0.333909, 0.452289) & (144, 144, 144) \\
(1/6, 1/3, 1/2) & \emph{No}  & (0.125793, 0.387562, 0.486645) & (88, 171, 155)\\
(1/2, 1/3, 1/6) & Yes & (0.390715, 0.371307, 0.237978) &(226, 151, 75)\\
\hline
\multicolumn{4}{|c|}{\bfseries QED Regime}\\
\hline
(1/3, 1/3, 1/3) & Yes & (0.214477, 0.332829, 0.452694) & (163, 163, 163) \\
(1/6, 1/3, 1/2) & \emph{No}  & (0.122010, 0.377990, 0.500000) & (98, 189, 180)\\
(1/2, 1/3, 1/6) & Yes & (0.388769, 0.371751, 0.239480) &(256, 171, 85)\\
\hline
\multicolumn{4}{|c|}{\bfseries QD Regime}\\
\hline
(1/3, 1/3, 1/3) & Yes & (0.225125, 0.335412, 0.439463) & (181, 181, 181) \\
(1/6, 1/3, 1/2) & \emph{No}  & (0.125416, 0.374584, 0.500000) & (107, 206, 205)\\
(1/2, 1/3, 1/6) & Yes & (0.398613, 0.369833, 0.231553) &(282, 188, 94)\\
\hline
\end{tabular}
\caption{Results from the numerical minimization of $\Delta$ for Example 1.}
\label{tbl:example1thetas}
\end{table}

\begin{figure}[htb]
\begin{minipage}{0.32\textwidth}\centering
$\theta_s=(1/3, 1/3, 1/3)$\\
\includegraphics[width=\linewidth]{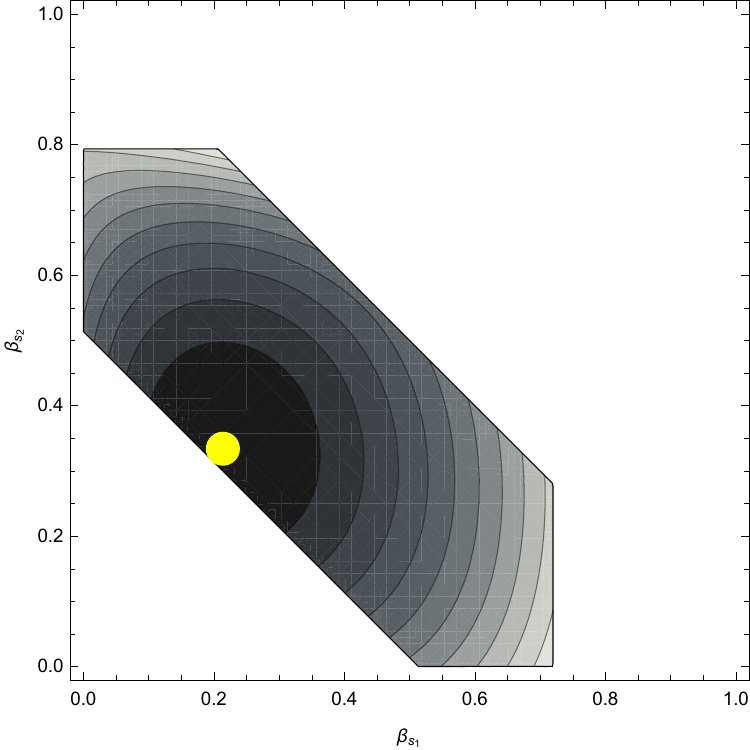} \\
(a)
\end{minipage}
\hfill\begin{minipage}{0.32\textwidth}\centering
$\theta_s=(1/6, 1/3, 1/2)$\\
\includegraphics[width=\linewidth]{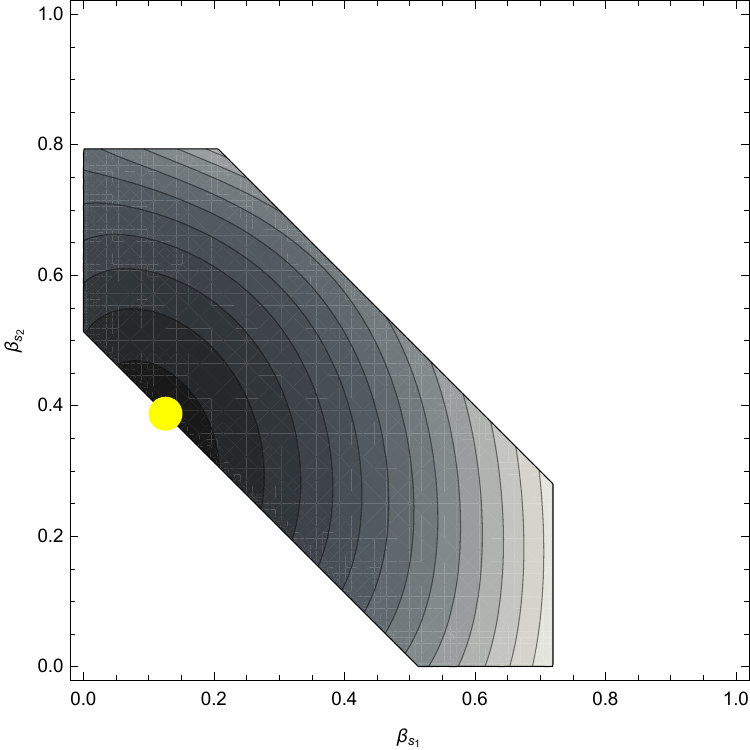} \\
(b)
\end{minipage}
\hfill\begin{minipage}{0.32\textwidth}\centering
$\theta_s=(1/2, 1/3, 1/6)$\\
\includegraphics[width=\linewidth]{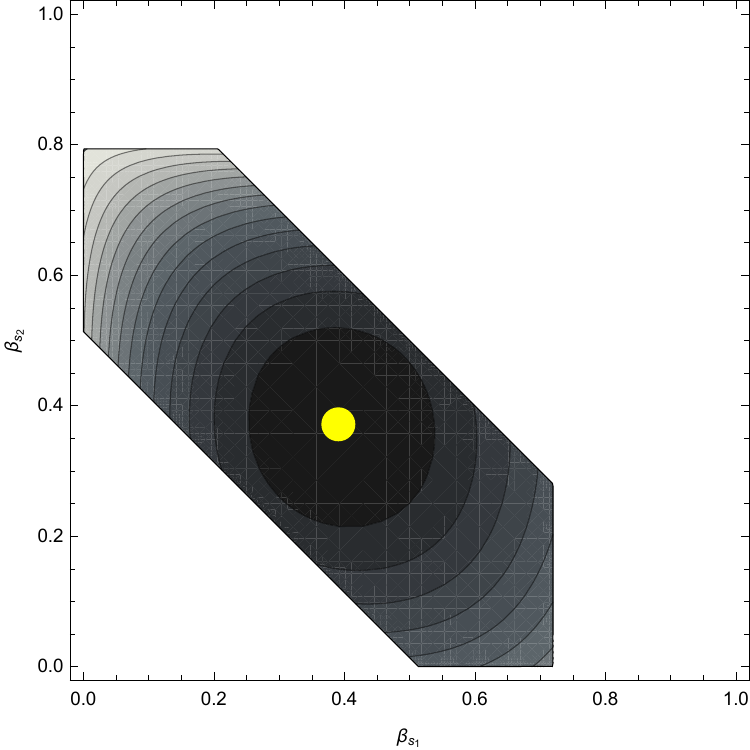} \\
(c)
\end{minipage}
\caption{Contour plots of $\Delta$ as a function of $\beta_{s_1}$ and $\beta_{s_2}$ for Example 1, in the ED regime. The highlighted area indicates the combinations of $\beta_{s_1}$ and $\beta_{s_2}$ that result in a system with complete resource pooling.}
\label{fig:example1thetasContour}
\end{figure}

Interestingly, the desired fractions $n_j/n$ for $j=1,2,3$ cannot be achieved for the case where $\theta_s=(1/6, 1/3, 1/2)$. In all three regimes, the global minimum of the function $\Delta$ is not attained within the region of complete resource pooling. This is illustrated in Figure \ref{fig:example1thetasContour}, where contour plots for $\Delta$ as a function of $\beta_{s_1}$ and $\beta_{s_2}$ are shown for the three different $\theta_s$ vectors, for the ED regime. The plots for the QD and QED regimes are omitted, as they look similar.

Since the global minimum of $\Delta$ for the desired vector $\theta_s=(1/6, 1/3, 1/2)$ is located outside the boundaries of the complete resource pooling region, we find that $\Delta>0$ for the vector $\beta_{s}$ that minimizes $\Delta$ within the complete resource pooling region. As a consequence, the relative staffing levels $n_j/n$ resulting from the vector $\beta_s$ differ from the desired $\theta_s$. For example, in the ED regime the recommended staffing levels are $(88, 171, 155)$, with relative values $(0.21349, 0.411097, 0.375413)$, whereas the desired $\theta_s=(1/6, 1/3, 1/2)$.


\paragraph{Example 3 revisited.} We now consider the $6 \times 6$ symmetric degree~3 graph from Example~3, with pooled service. As an illustration, we distinguish between three different combinations of $\theta_{s_j}$'s, similar to the previous example:
\begin{itemize}
\item $\theta_{s_1}=\theta_{s_2}=\dots=\theta_{s_6}=1/6$, i.e. all equal,
\item $\theta_{s}=(1,2,3,4,5,6)/21$, i.e. increasing in the server type,
\item $\theta_s=(6,5,4,3,2,1)/21$, i.e. decreasing in the server type.
\end{itemize}
Using the same numerical optimization routine, we determine the required $\beta_{s_j}$'s from which the staffing levels $n_{s_j}$ can be determined for the QD, QED, and the ED regime. The results can be found in Table \ref{tbl:example3thetas}. Some interesting conclusions can be drawn from this table.
First, we see that all three designs are realizable with complete resource pooling.  This may be explained by the added amount of flexibility allowed by having each server compatible with 3 types of customers, i.e. large degree of overlap in compatibilities.
Second, it is immediately clear that the relative staffing levels $n_j/n$ resulting from the numerical procedure are all equal to the specified values $\theta_{s_j}$, meaning that the minimization of $\Delta$ resulted in a global minimum of 0 satisfying the conditions for complete resource pooling.
Third, it can be seen that the vector $\beta_s$ that minimizes $\Delta$ is the same for all three regimes. This can easily be explained from the fact that all patience distributions are the same, meaning that the values of $\alpha_{c_i}$, $i=1,\dots,I$, do not depend on the selected regime and neither does the function $\Delta$.

\begin{table}[htb]
\begin{tabular}{|c|c|c|c|}
\hline
Desired $\theta_{s}$ & Realized & $\beta_s$ & Staffing levels for $\lambda=100$ \\
\hline
\multicolumn{4}{|c|}{\bfseries ED Regime}\\
\hline
(1,1,1,1,1,1)/6  & Yes & (0.147, 0.187, 0.147, 0.187, 0.147, 0.187) & (52, 52, 52, 52, 52, 52)\\
(1,2,3,4,5,6)/21 & Yes & (0.041, 0.131, 0.117, 0.198, 0.200, 0.312) & (15, 29, 44, 59, 74, 88)\\
(6,5,4,3,2,1)/21 & Yes & (0.251, 0.264, 0.175, 0.165, 0.094, 0.052) & (92, 77, 61, 46, 31, 15)\\
\hline
\multicolumn{4}{|c|}{\bfseries QED Regime}\\
\hline
(1,1,1,1,1,1)/6  & Yes & (0.147, 0.187, 0.147, 0.187, 0.147, 0.187)  & (57, 57, 57, 57, 57, 57)\\
(1,2,3,4,5,6)/21 & Yes &(0.041, 0.131, 0.117, 0.198, 0.200, 0.312) & (16, 33, 49, 65, 81, 98)\\
(6,5,4,3,2,1)/21 & Yes & (0.251, 0.264, 0.175, 0.165, 0.094, 0.052) & (102, 85, 68, 51, 34, 17)\\
\hline
\multicolumn{4}{|c|}{\bfseries QD Regime}\\
\hline
(1,1,1,1,1,1)/6  & Yes & (0.149, 0.184, 0.149, 0.184, 0.149, 0.184)  & (66, 66, 66, 66, 66, 66)\\
(1,2,3,4,5,6)/21 & Yes & (0.042, 0.126, 0.120, 0.197, 0.205, 0.309) & (19, 37, 56, 75, 94, 112)\\
(6,5,4,3,2,1)/21 & Yes & (0.255, 0.260, 0.177, 0.162, 0.094, 0.052) & (116, 97, 78, 58, 39, 19)\\
\hline
\end{tabular}
\caption{Results from the numerical minimization of $\Delta$ for Example 3.}
\label{tbl:example3thetas}
\end{table}

\section{Discussion}
\label{sec.discussion}
Our purpose in this paper was twofold:  To verify a conjecture on matching rates, and to provide a useful tool for design of parallel service systems.
We now discuss these two points.

Our computational results seem to indicate that our conjecture on matching rates may be correct.  In fact our simulations show much more than that.  Our conjecture is that under many server scaling matching rates will converge precisely to those of
the infinite matching model.
This would in particular imply that our designs would converge to a deterministic limit, in which for ED designs all patient customers wait exactly $W$, in QD design all servers will idle exactly $T$, and in QED nobody would idle or wait.  Our simulations indicate that this may be true, and leave the question of proof of the limiting result open.

However, the simulations also indicate that this convergence does not require unreasonably large $n$.  In fact, for quite realistic values of $\lambda$ and $n$ we find that our designs perform extremely well:  The matching rates in the simulations are very close to those predicted by the approximations, but more important, the performance of the systems is very close to the required service quality and utilization parameters as specified in the designs.

We have found that in  ED mode the fraction of abandonments is almost precisely the pre-specified value, and that the waiting time distribution of patient customers is distributed around the mean value $W$ which is the conjectured limiting value for $n\to\infty$.   We have also shown that in QD mode the idle times are distributed around the mean value $T$ which is the conjectured limiting value for $n\to\infty$.
We have also shown that in the QED mode the system works in perfect balance between customers and servers, with complete resource pooling and uniform service level for all customers, with negligible abandonments.  All this for the whole range of values of $\lambda$, from 20 to 200.

We note again that our conjecture on matching rates is for general bipartite compatibility graphs, general patience and service time distributions, and general renewal arrivals.
Our simulations are of course limited and cannot replace a mathematical proof, and we see no way of making them extensive enough to validate our designs in all situations.    However the simulations do cover various graphs, a wide range of distributions, and  realistic quality parameter values, and as can be seen, the results agree with our predictions in all the examples we tried.

It is important to reiterate that these systems under FCFS-ALIS are completely intractable.  without recourse to our matching rate calculations it would be impossible to design a system that would achieve resource pooling and approximate our design parameters.  The paper of Foss and Chernova \cite{foss-chernova:98} abundantly illustrates this intractability.

Extensive simulation confirmed that the algorithms based on our conjecture are accurate and effective: they produce work force levels and, in case of differentiated service, a redesigned compatibility graph that meet targeted quality of service requirements.
Moreover, the heuristic algorithms also appeared to work well when the required work force levels are not so large.
As such, these algorithms provide a valuable tool to support decisions on the design of multi-type parallel service systems.

We realize that in real large scale systems there are many specific needs and constraints that will make pure FCFS-ALIS policy impossible to use.  However, we believe that even if FCFS-ALIS is not used throughout, it is used for a significant  fraction of the scheduling decisions in real large systems.  Currently these systems are often evaluated and redesigned based on simulation studies.  We suggest that our tools that allow reasonably accurate evaluation of system performance under FCFS-ALIS using direct computation, can extend the range of design tools based on simulations, and provide a useful method of design.
 In particular, our tools can provide quick answers on how system performance will change if we change the compatibility graph, if we redefine service priorities, if we change quality of service parameters, if we reallocate server types, or if we change service rates by redefining tasks.

Much still remains to be done in this stream of research.  A proof of the conjecture,  in general form or under more limiting assumptions, is the first task.  A preliminary attempt to verify the conjecture for many server N-system in the Poisson-exponential case is presented in \cite{zhan-weiss:16}.  It would also be interesting to  evaluate the efficacy of FCFS-ALIS policy in comparison with other policies.  Finally,  application to a real system with real data will be most illuminating.


\end{document}